\documentclass[reqno,11pt]{amsart}
\usepackage{amsmath,amsthm,amscd,amsfonts,amssymb,hyperref,floatrow}
\usepackage{graphicx, color,dsfont,xcolor}

\numberwithin{equation}{section}

\newtheorem{theorem}{Theorem}[section]
\newtheorem{lemma}[theorem]{Lemma}
\newtheorem{proposition}[theorem]{Proposition}
\newtheorem{corollary}[theorem]{Corollary}
\newtheorem{rem}[theorem]{Remark}

\DeclareMathSymbol{\leqslant}{\mathalpha}{AMSa}{"36} 
\DeclareMathSymbol{\geqslant}{\mathalpha}{AMSa}{"3E} 
\DeclareMathSymbol{\eset}{\mathalpha}{AMSb}{"3F}     
\renewcommand{\leq}{\;\leqslant\;}                   
\renewcommand{\geq}{\;\geqslant\;}                   

\newcommand{\R}{\mathbb{R}}

\newcommand{\N}{\mathbb{N}}

\def \P{ \mathbb P  }
 \newcommand{\eps}{\varepsilon}
\def \E{ \mathbb E  }

\begin{document}
\title{Tracy-Widom at high temperature}

\author{Romain Allez \and Laure Dumaz}
\address{Weierstrass Institute, Mohrenstr. 39, 10117 Berlin, Germany.}
\address{Statistical Laboratory, Centre for Mathematical Sciences, Wilberforce Road, Cambridge, CB3 0WB, United Kingdom. }

\email{romain.allez@gmail.com,L.dumaz@statslab.cam.ac.uk} 
\date{\today}

\maketitle

\begin{abstract}
We investigate the marginal distribution of the bottom eigenvalues of the stochastic Airy operator 
when the inverse temperature $\beta$ tends to $0$. 
We prove that the minimal eigenvalue, whose fluctuations are governed by the Tracy-Widom $\beta$ law, 
converges weakly, when properly centered and scaled, to the Gumbel distribution. 
More generally we obtain the convergence in law of the marginal distribution 
of any eigenvalue with given index $k$.  
Those convergences are obtained after a careful analysis of the explosion times process 
of the Riccati diffusion associated to the stochastic Airy operator. 
We show that the empirical measure of the explosion times converges 
weakly to a Poisson point process using estimates proved in [L. Dumaz and B. Vir\'ag.  Ann. Inst. H. Poincar\'e Probab. Statist. {\bf 49}, 4, 915-933, (2013)].  
We further compute the empirical eigenvalue density of the stochastic Airy ensemble on the macroscopic scale when $\beta\to 0$. 
As an application, we investigate the maximal eigenvalues statistics 
of $\beta_N$-ensembles when the repulsion parameter $\beta_N\to 0$ when $N\to +\infty$. 
We study the double scaling limit $N\to +\infty, \beta_N \to 0$ and argue with heuristic and numerical arguments 
that the statistics of the marginal distributions can be deduced following the ideas of 
[A. Edelman and B. D. Sutton. J. Stat. Phys. {\bf 127} 6, 1121-1165 (2007)] and [J. A. Ram\'irez, B. Rider and B. Vir\'ag. J. Amer. Math. Soc. {\bf 24} 919-944 (2011)] from our later study of the stochastic Airy operator.

\end{abstract}
\vspace{1cm}

\section{Introduction}
One of the most influential random matrix theory (RMT)
developments of the last decade was discovery in 2002 of the so-called tridiagonal 
$\beta$-ensembles by Dumitriu and Edelman in \cite{dumitriu}. 
The tridiagonal random matrices of this ensemble 
have explicit and independent (up to symmetry) entries and their eigenvalues are distributed 
according to the equilibrium joint probability density function (jpdf) of charged particles in a 
one dimensional Coulomb gas with electrostatic repulsion, confined in a quadratic potential and subject to a thermal 
noise at temperature $T=1/\beta$ for {\it arbitrary} $\beta>0$. More precisely, the jpdf of the eigenvalues is given by
\begin{equation}\label{eq.P_beta}
P_\beta(\lambda_1,\cdots,\lambda_N) = \frac{1}{Z_N^\beta} \prod_{i < j} |\lambda_i-\lambda_j|^{\beta} \exp(-\frac{1}{4}\sum_{i=1}^N \lambda_i^2)\,. 
\end{equation}
For $\beta=1$ (respectively $\beta=2,4$), this jpdf arises as the joint law of the eigenvalues of the classical Gaussian orthogonal 
(respectively unitary, symplectic) ensembles, whose linear eigenvalues statistics were extensively studied in 
the literature (see \cite{agz,silverstein,mehta,forrester,handbook} for a review of RMT and its applications). 

The introduction of the tridiagonal random matrices for {\it arbitrary} $\beta$-ensembles 
has led to considerable progress for the study of linear statistics of the point process with jpdf $P_\beta$. 
They have permitted to prove that the largest eigenvalues converge jointly in distribution to the low-lying eigenvalues of the random Schrodinger
operator, also called {\it stochastic Airy operator}, $-\frac{d^2}{dt^2}+ t +\frac{2}{\sqrt{\beta}} b'_t$, restricted to the positive half-line, 
where $b'$ is a white noise on $\R_+$ \cite{edelman,virag} ( see also \cite{schehr} for a review on the top eigenvalue statistics of random matrices). 

In this paper, we are interested in the limiting marginal distributions of the bottom eigenvalues of the stochastic Airy operator
when the parameter $\beta$ tends to $0$.  
As mentioned above, the stochastic Airy operator appears as the continuous scaling limit of $\beta$-ensembles at the edge of the spectrum 
and we shall see that the question we investigate is in fact related to the 
largest eigenvalues statistics of $\beta$-ensembles 
when the $\beta$ index scales with the dimension $N$ such that $\beta_N \to 0$ when $N\to +\infty$.

A somehow related question was investigated in \cite{jp-alice} (see also \cite{alice}) where the authors
consider the empirical eigenvalues density of $\beta$-ensembles in the limit of large dimension $N$ and with $\beta=2c/N$. 
The limiting family of probability density $\{\rho_c, c\geq 0\}$ is computed explicitly in terms of Parabolic cylinder functions and is proved 
to interpolate continuously between the Gaussian shape (obtained for $c=0$) 
and the Wigner semicircle shape (which is recovered 
when $c\rightarrow +\infty$). The $\beta$-Wishart ensemble was handled similarly in \cite{satya}.

The question of the characterization of an interpolation between the Tracy-Widom $\beta$ distribution (which governs the typical 
fluctuations of the 
top eigenvalue as $N \to \infty$ with $\beta>0$ fixed)
and the Gumbel distribution (which governs the typical fluctuations 
of the maximum of independent Gaussian variables -- corresponding to the $\beta=0$ case) was raised in \cite{jp-alice} and \cite{johansson}.

We answer this question proving that the Tracy-Widom $\beta$ distribution 
converges weakly (when properly rescaled and centered) to the Gumbel 
distribution when $\beta$ goes to $0$. This is the content of Theorem \ref{TracyWidom.Gumbel}.  
We use the characterization of the marginal distributions of the low lying eigenvalues 
of the stochastic Airy operator in terms of the explosion times process of the associated Riccati diffusion \cite{virag}. 
We show that the empirical measure of the explosion times converges weakly in the space of Radon measures 
to an inhomogeneous Poisson point process on $\R_+$ with explicit intensity. 
The weak convergence 
of all the marginal distributions of the second, third, etc eigenvalues can be readily deduced. 
Although we expect the minimal eigenvalues to have Poissonian statistics in the small $\beta$ limit, 
the convergences of the joint distribution of the $k$ bottom eigenvalues for any fixed index $k$  
seem to be difficult to prove as there is not a simple characterization of this law in terms of a single diffusion.  
It is still characterized in this setting in terms of a family of coupled diffusions but the interaction between those diffusions is complex and makes the analysis difficult (see Figures \ref{explosion1} and \ref{couplage_diffusions} below).

As an application, we investigate (with heuristic and numerical arguments)
the weak convergence of the top eigenvalue of $\beta_N$-ensembles
in the double scaling limit $N\to +\infty$ and $\beta_N\to 0$. 
We revisit the ideas of \cite{edelman} which proposes that tridiagonal random matrices of $\beta$-ensembles 
are properly viewed as finite difference schemes of the stochastic Airy operator. 
From our heuristic discussion in section \ref{rmt}, this relation seems to remain valid also in the regime $\beta_N\to 0$ and 
permits to establish the weak convergence of the top eigenvalue
of $\beta_N$ ensembles to the Gumbel distribution, for any sequence $\beta_N$ such that 
\begin{equation*}
1\, \gg \, \beta_N \,  \gtrapprox  \, \frac{\ln N}{N} \,.
\end{equation*}
We explicit the scaling and centering of this convergence which are in fact the same as in the convergence of the Tracy-Widom 
$\beta$ distribution to the Gumbel distribution. 
Again, the allied results on the second, third, fourth, etc eigenvalues can also be derived from our former results on the Stochastic Airy ensemble
(${SAE}_\beta$ for short). 

We finally mention that our derivation does not cover the case where $\beta\sim 1/N$, 
which is highly interesting as the interpolation for the empirical spectral distribution occurs on this range of $\beta_N$ \cite{jp-alice}. 
For $\beta_N$ decreasing as slowly as $\ln N/N$ or even more slowly, the typical fluctuations of the top eigenvalues seem to 
enter the Gumbel regime in the sense that the centerings and scalings are found to be the same as in the classical setting of independent 
Gaussian variables.

\medskip

To facilitate the reading, let us draw up a short outline of the paper. 
In Section \ref{sae.review}, we give a brief review on Stochastic Airy Ensemble (${\rm SAE}_\beta$ for short), 
recalling in particular the correspondence 
between the law of the eigenvalues and the law 
of the explosion times process of the associated Riccati diffusion established by \cite{virag}. 

In section \ref{exit.time.statio}, we revisit the classical problem of the exit time from a domain 
of a diffusion which evolves in a {\it stationary} potential. This problem 
is the stationary counterpart of our main study and turns to be useful in the next sections 
to approximate the non stationary Riccati diffusion and in particular its explosion times. 
We provide a simple characterization of the law of the exit time which permits to prove its weak convergence 
to an exponential distribution,
when the trap gets very deep in comparison to the noise. Then we consider the explosion times process of the stationary Riccati diffusion 
and we prove that 
it converges to a (homogeneous) Poisson point process. Finally, we discuss in view of subsection \ref{macro.density}  
the Fokker Planck equation which relates 
between the transition probability distribution of a diffusion and the flux of probability in the system.

In Section \ref{sae.results}, we state our results on the 
convergence of the distribution of the minimal eigenvalue of the stochastic Airy operator, i.e. of the Tracy-Widom $\beta$ law, 
to the Gumbel law.
This is straightforwardly deduced from the convergence of the explosion time process of the diffusion. 
The convergences of the marginal distributions of the other neighboring minimal eigenvalues can be deduced as well. 
The proofs of those results appear in subsection \ref{main.proof}. At the end of the section, we compute with a perturbative \emph{heuristic}
method the empirical eigenvalue density of the stochastic Airy operator as $\beta \to 0$ on the macroscopic scale, i.e. 
without any zooming in the minimal eigenvalues scaling region. 

As an application of our results, we discuss in 
section \ref{rmt} the marginal statistics of the 
minimal eigenvalues of $\beta_N$-ensembles in the double scaling limit $\beta_N\to 0, N \to +\infty$. 
We conjecture that they can be readily deduced from our results since the tridiagonal random matrices, when zooming 
in the edge scaling region, are well approximated, even when $\beta_N\to 0$, 
by the stochastic Airy operator. 
Some technical computations and 
auxiliary proofs are gathered into appendices. 

\vspace{0.5cm}

{\bf Acknowledgments.}
We are very grateful to J.-P. Bouchaud for many enlightening discussions. 
We have shared many ideas with him and benefited from his insight along this project. 

We also thank Edouard Br\'ezin for his explanations on \cite{bb},  Antoine Levitt and Alain Comtet for useful 
comments on Schrodinger operators and for pointing out reference \cite{texier} as well as
Satya N. Majumdar for interesting discussions related to problems at stake in section \ref{rmt}. 

RA received funding from the European Research Council under the European
Union's Seventh Framework Programme (FP7/2007-2013) / ERC grant agreement nr. 258237 and thanks the Statslab in DPMMS, Cambridge for its hospitality at the time this work was finished.
The work of L.D. was supported by the 
Engineering and Physical Sciences
Research Council under grant EP/103372X/1 and L.D. thanks the hospitality of the maths department of TU and the Weierstrass institute in Berlin.

\section{Stochastic Airy operator: A short review} \label{sae.review}
In the first subsection, we recall the recent results on the spectral statistics of the stochastic Airy operator 
obtained in \cite{edelman,virag,laure,celine,forrester.tw}.  
Then, we give a characterization of the first exit time (also called blow-up time) of the non homogeneous {\it Riccati} diffusion 
associated to the stochastic Airy operator.

The stochastic Airy operator $\mathcal{H}_\beta$ is defined formally (see \cite[Section 2]{virag} for a precise definition) for $\beta>0$ and $t\geq 0$ as 
\begin{equation}\label{def.Abeta}
\mathcal{H}_\beta := -\frac{d^2}{dt^2} + t + \frac{2}{\sqrt{\beta}} \, B'_t
\end{equation}
where $B'_t$ is a white noise on $\R_+$.  Following \cite{virag}, denote by $S^*$ the space of functions $f$ satisfying $f(0)=0$
and $\int_0^{\infty} (f')^2 +(1+t) f^2 < \infty$. 
We will say that $(\phi_\lambda,\lambda)\in S^*\times \R$ is an eigenfunction/ eigenvalue pair for $\mathcal{H}_\beta$ if $||\phi_\lambda||_2=1$ 
and if
\begin{equation}\label{eq.diff.phi}
\phi_\lambda''(t) = (t-\lambda) \phi_\lambda(t) + \frac{2}{\sqrt{\beta}}\, \phi_\lambda(t) \, B'_t 
\end{equation}
holds for all $t\geq 0$ in the following integration by part sense, 
\begin{equation}\label{eq.diff.phi.2}
\phi_\lambda'(t) - \phi_\lambda'(0) = \int_0^t (s-\lambda)\, \phi_\lambda(s) \, ds + \frac{2}{\sqrt{\beta}} \left(B_t \, \phi_\lambda(t) - \int_0^t B_s \,  \phi_\lambda'(s) \, ds \right)\,. 
\end{equation}

Ram\'irez, Rider and Vir\'ag proved in \cite[Theorem 1.1]{virag} that, almost surely, for each $k\geq 0$, the set of eigenvalues of $\mathcal{H}_\beta$ 
has a well defined $(k+1)$ st lowest element denoted $\Lambda_k^\beta$.  
Furthermore the law of any eigenvalue $\Lambda_k^\beta$ with given index $k$
is characterized in term of the explosion times of a stochastic process $(X_\lambda(t))_{t\geq 0}$ 
defined through the Riccati change of functions $X_\lambda(t) := \phi_\lambda'(t)/\phi_\lambda(t)$.
This stochastic process is a diffusion process
whose initial condition and Langevin equation are 
obtained from the Dirichlet boundary condition $\phi_\lambda(0)=0$ \footnote{Necessarily $\phi_\lambda'(0)\neq 0$ (otherwise $\phi_\lambda$ is identically $0$) 
and the two signs of $\phi_\lambda'(0),\phi_\lambda(0+)$ are equal.} and \eqref{eq.diff.phi} 
\begin{equation}\label{eq.lang.xt}
X_\lambda(0) = +\infty \quad {\rm and} 
\quad dX_\lambda(t) = \left(t - \lambda - X_\lambda(t)^2\right) \, dt \, + \frac{2}{\sqrt{\beta}}\, dB_t  \quad {\rm for}\quad  t\geq 0,
\end{equation}
where $B_t$ is a standard Brownian motion. 
Solutions of \eqref{eq.lang.xt} may blow up to $-\infty$ at finite times, as will happen whenever 
$\phi_\lambda$ vanishes. In this case, the diffusion $X_\lambda$ immediately restarts at $+\infty$ at that time in order to continue the solution 
corresponding to the Langevin equation \eqref{eq.lang.xt}.

The authors of \cite{virag} 
prove that the operator $\mathcal{H}_\beta$
satisfies a Sturm-Liouville like property in the sense that {\it the number of eigenvalues of $\mathcal{H}_\beta$ at most $\lambda$ 
is equal to the total number of explosions of the diffusion $(X_\lambda(t))$ on $\R_+$} \footnote{It is also proved 
that the diffusion $X_\lambda$ has only a finite number of explosions in $\R_+$ and thus that the number of eigenvalues
below $\lambda$ is almost surely finite.}.  

Before explaining the idea behind this key relation, let us preliminary state the so-called {\it increasing property} of the 
coupled family of diffusions $(X_\lambda)_{\lambda\in \R}$, which 
will be used many times in the paper. The increasing property is rather intuitive and can be enunciated as follows: 
if $\lambda' \leq \lambda$, then the number of explosions of $X_{\lambda}$ is stochastically bounded above by the number of explosions 
of $X_{\lambda'}$ on any compact interval $[0,T]$. Equivalently, the diffusion $X_{\lambda'}$ remains below $X_{\lambda}$ until its 
first explosion time and can not cross the trajectory of $X_\lambda$ from below to above. Intuitively
this property is rather obvious since the drift of $X_{\lambda'}$ pulls stronger downside than the one of $X_{\lambda}$. 
Note however that the comparison theorem for sdes (see \cite[Proposition 2.18]{karatzas} or \cite[Theorem (3.7), Chapter IX]{yor}) does not apply directly (the drifts are not Lipschitz) and one needs to use a localization argument
before applying it. 

\begin{figure}[h!btp] 
     \center
     \includegraphics[scale=0.7]{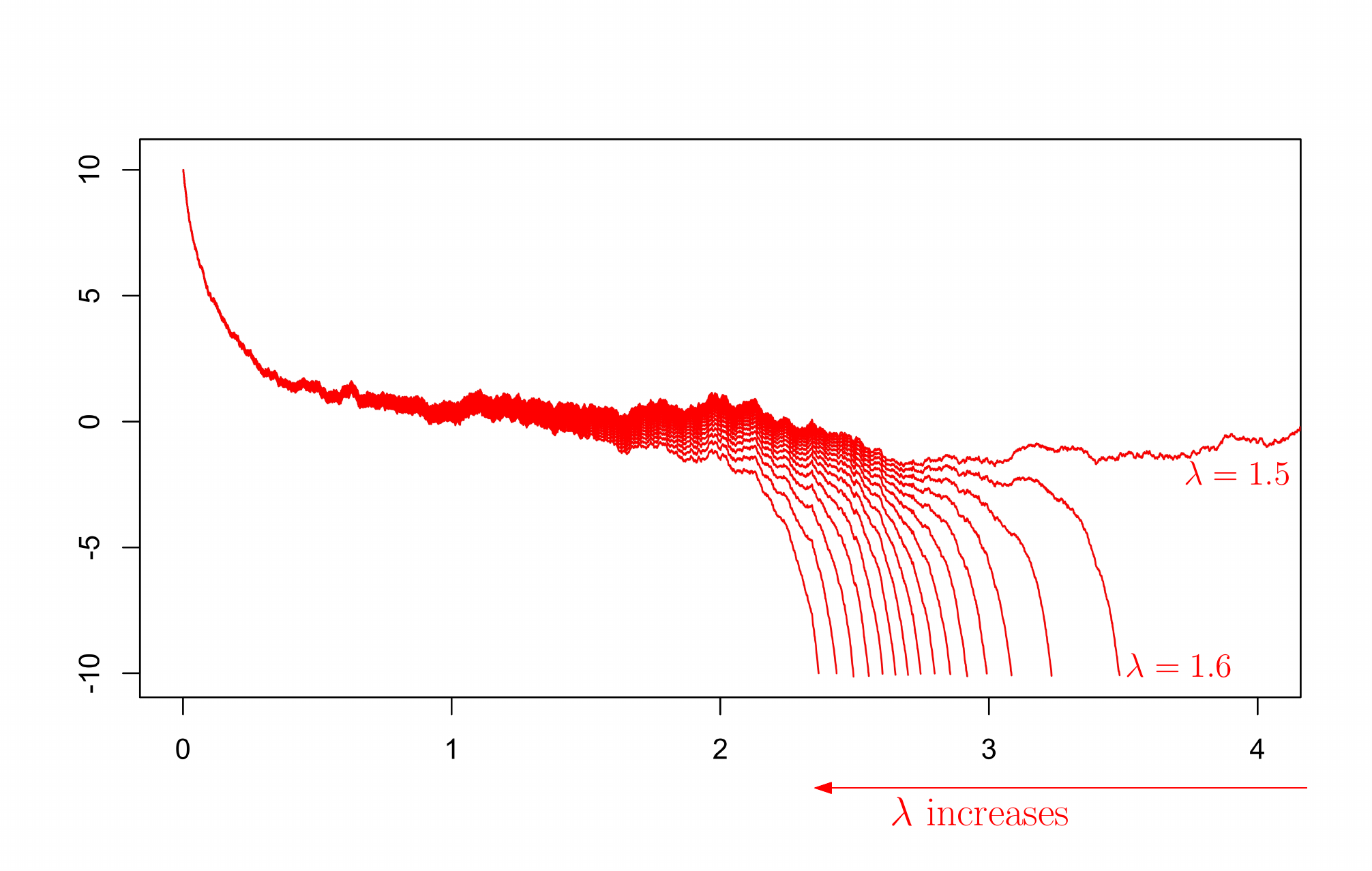}
\caption{Simulated paths of diffusions $X_\lambda$ driven by the same Brownian motion for several values of $\lambda$. 
The values of $\lambda$ are between $1.5$ and $3$ on a grid of mesh $10^{-1}$. We took $\beta=4$. 
The smallest eigenvalue $-TW(4)$ of the Airy operator is between $1.5$ (no explosion) and $1.6$ (at least one explosion) on this event.}
\label{explosion1}
\end{figure}

\begin{figure}[h!btp] 
     \center
     \includegraphics[scale=0.6]{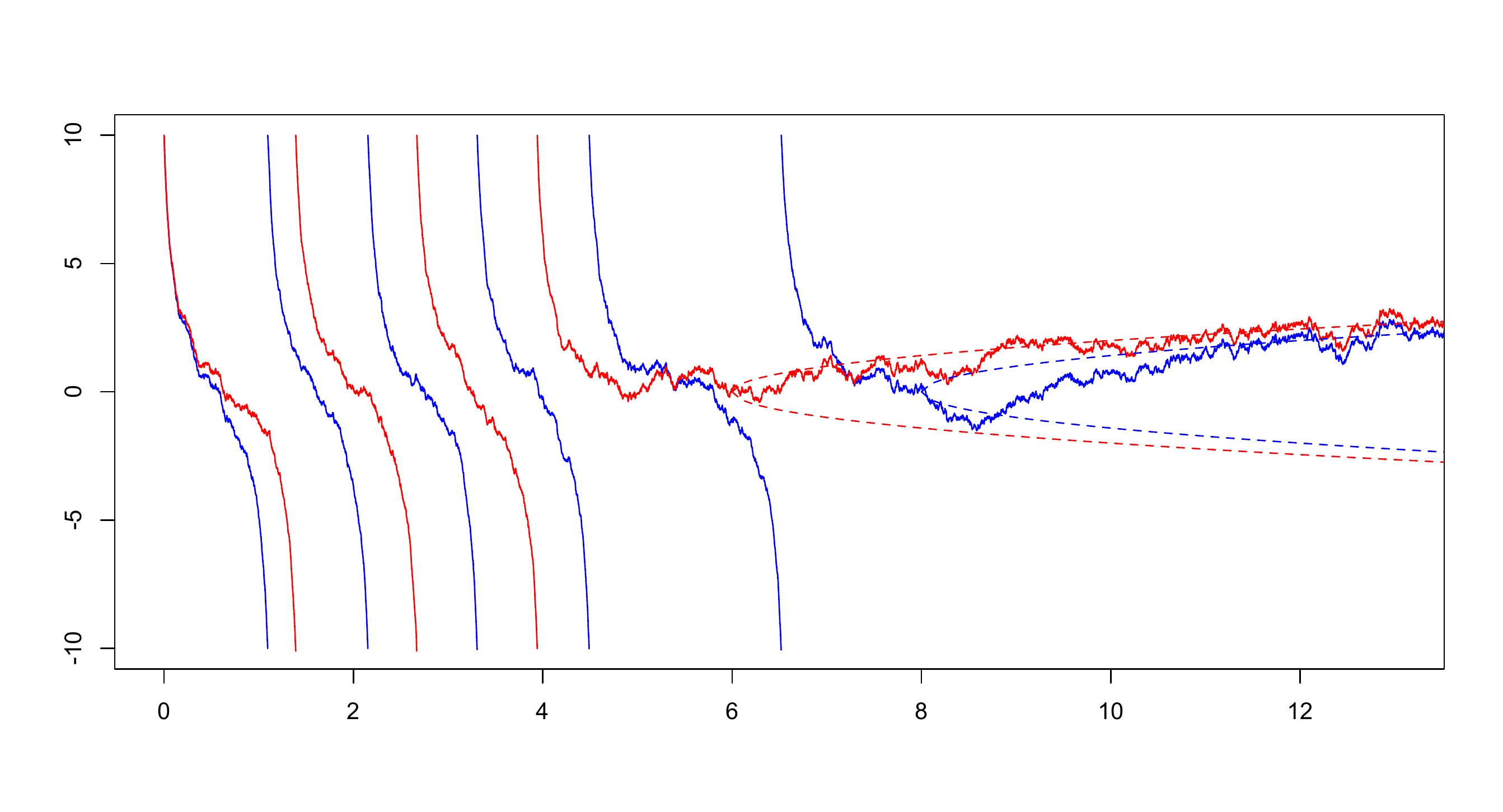}
\caption{(Color online) Simulated paths of two diffusions $X_\lambda$ and $X_{\lambda'}$ driven by the same Brownian motion 
with $\lambda=8$ (blue curve) and $\lambda'= 6$ (red curve), for $\beta=4$. We add in dashed lines the corresponding 
parabolas where the drifts cancel. When $t\to +\infty$, the diffusions converge to the upper part of their respective parabolas. 
On this event, we have $\Lambda_4^\beta < 8< \Lambda_5^\beta$ and $\Lambda_2^\beta < 6<\Lambda_3^\beta$.}
\label{couplage_diffusions}
\end{figure}

Now we explain the key relation. First, we will look at the operator $\mathcal{H}_\beta^L$ defined on the set of functions with support 
in the truncated interval $[0, L]$ with Dirichlet boundary conditions at both endpoints. 
This truncated operator is shown to approximate closely
the operator $\mathcal{H}_\beta$ when $L \to \infty$ in a precise topology \cite{virag}.
At fixed $L$, the definitions of the eigenfunctions and their associated Riccati transform imply that
$\lambda\in \R$ is an eigenvalue of the operator
$\mathcal{H}_\beta^L$ \emph{if and only if} the Riccati diffusion $X_{\lambda}$ explodes precisely at the end 
point $L$. We can deduce that 
the number of eigenvalues of the operator $\mathcal{H}_\beta^L$ smaller than $\lambda$ equals the number of real numbers 
$\lambda'\le \lambda$ such that $(X_{\lambda'})$ blows up exactly in $L$. Indeed, if $\lambda'$ slowly decreases 
from $\lambda$ to $-\infty$, we can see from the {\it increasing property} 
and the continuity of the whole trajectory $(X_{\lambda'}(t))$ with respect to 
$\lambda'$ that the explosion waiting times of $X_{\lambda'}$ will simultaneously get longer and longer  
to finally diverge to $+\infty$ when $\lambda'$ gets too small for the diffusion to explode in $\R_+$ 
(see Fig. \ref{explosion1} and \ref{couplage_diffusions}). 
In particular, any explosion of $X_{\lambda}$ which occurs before time $L$ will slowly translate continuously to occur at some point 
exactly in $L$ 
for some $X_{\lambda'}$ with $\lambda'\leq \lambda$. Thus the number of eigenvalues of $\mathcal{H}_\beta^L$ smaller than $\lambda$ 
is equal to the number of explosions of $X_\lambda$ on the interval $[0;L]$.  
Taking $L \to \infty$ finally establishes this property for the operator $\mathcal{H}_\beta$.


The marginal laws of the eigenvalues are 
characterized in a rather simple way in 
term of the distribution of the explosion times process of one single Riccati diffusion $X_\lambda$ as a function of the parameter 
$\lambda$. We shall use this characterization several times throughout the paper. 

In particular,
the cumulative distribution of the lowest eigenvalue $\Lambda_0^\beta$ satisfies  
\begin{equation}\label{charac.tw}
\P\left[\Lambda_0^\beta < \lambda \right] = \P\left[ X_\lambda(t) \,  \textrm{blows up to $- \infty$ in a finite time} \right] \,.
\end{equation}
And for $k\geq 1$, we have 
\begin{align}
\P\left[\Lambda_k^\beta < \lambda \right] = \P\left[ X_\lambda(t) \,  \textrm{blows up to $- \infty$ at least $k+1$ times} \right] \,.
\end{align}

The joint distribution of the $k$ bottom eigenvalues can also be characterized in this setting but 
in a more complicated way from the law of the family 
of coupled diffusions $((X_\lambda(t))_{t \geq 0},\lambda\in \R)$ all satisfying \eqref{eq.lang.xt} with the {\it same} 
driving Brownian motion $(B(t))_{t\geq 0}$. 
For instance, with $\lambda' < \lambda$,
\begin{align}\label{coupling}
\P&[\Lambda_{k-1}^\beta \leq \lambda, \Lambda_{k}^\beta > \lambda' ] \notag\\ 
&= \P\bigg[ X_{\lambda} \mbox{ blows up to $- \infty$ at least $k$ times }, X_{\lambda'} \mbox{ blows up to $- \infty$ at most $k$ times }\bigg]\,. 
\end{align}  
Although the two diffusions $X_\lambda$ and $X_{\lambda'}$ are driven by the same Brownian motion, they have different drifts
and (to our knowledge) this makes the right hand side probability of \eqref{coupling} difficult to estimate in practice.

In view of \cite[Theorem 1.1]{virag} which states that the top eigenvalue of the tridiagonal matrices 
properly centered and scaled converges in law to $-\Lambda_0^\beta$, the Tracy-Widom $\beta$ law ($TW(\beta)$ for short) 
has been defined for general $\beta>0$ as
the law of the random variable $-\Lambda_0^\beta$. 

Using the characterization \eqref{charac.tw} in terms of the diffusion $X_\lambda$, the first two leading terms of the right 
large deviation tail of the Tracy-Widom $\beta$ distribution were rigorously obtained by the second author and Vir\'ag in \cite{laure}. 
Those two first terms were also computed by Forrester in \cite{forrester.tw} using a different method.   
Finally, the right tail of the Tracy-Widom $\beta$ law was computed to all orders by Borot and Nadal in \cite{celine}
using heuristic arguments. 
Their result valid in the  $\lambda\rightarrow +\infty$ limit reads 
\begin{equation}\label{right.tail.tw}
\P\left[TW(\beta)> \lambda\right] = \frac{\Gamma\left(\frac{\beta}{2}\right)}{(4\beta)^{\frac{\beta}{2} } 2 \pi} \lambda^{-\frac{3\beta}{4}} 
\exp\left( -\frac{2}{3} \,\beta \, \lambda^{3/2} \right)  \exp\left( \sum_{m=1}^{+\infty} \frac{\beta}{2} R_m(\frac{2}{\beta}) \lambda^{-3m/2}  \right),
\end{equation} 
where the $R_m$ are (explicit) polynomials of degree at most $m+1$. 

It is in fact possible to obtain an analytical characterization of the law of the minimal eigenvalue $\Lambda_0^\beta$. 
In \cite[Theorem 1.7]{bloemendal}, the authors prove that the cumulative distribution 
of $\Lambda_0^\beta$ is the unique solution of a boundary value problem. 
Using similar techniques, we can prove slightly more: The \emph{Laplace transform} of the first explosion time is the solution of a boundary value problem as well. See Remark \ref{rem.LT} for more details. But it turns out that this boundary value problem is hard to analyse, even in the limit $\beta \to 0$, and this approach in spite of its robustness did not permit us to prove Theorem \ref{th.conv.explosion.times}.

\section{Trapping of a diffusion in a stationary well}\label{exit.time.statio}
We now revisit the classical problem of the exit time from a domain 
of a diffusion which evolves in a stationary potential. 
The {\it small noise limit} was widely studied in the literature (see e.g. \cite{freidlinwentzell, zeitouni, varadhan}) 
using large deviation theory. The results we derive in this section are not new and hold in a general setting. 

The diffusion $Y_a$ considered below in the small noise limit appears in the study \cite{texier,lloyd,halperin,mckean} 
of the law of the ground state (minimal eigenvalue) 
of the {\it Hill}'s operator, defined as
\begin{align*}
\mathcal{G}_L:= -\frac{d^2}{dt^2} + B'(t)
\end{align*}
where $B'(t)$ is as before a white noise on the segment $[0,L], L>0$. 
In this context, due to the stationarity (absence of the linear $t$ term), we need to restrict to a finite perimeter 
$L>0$ and we work with Dirichlet boundary conditions $\psi_a(0)=\psi_a(L)=0$ for the eigenvectors such that 
$\mathcal{G}_L \psi_a=a\psi_a$. 
As in the previous section, the law of the minimal eigenvalues 
$A_k^L$ of the operator $ \mathcal{G}_L$ can be characterized in terms of the family of diffusions  
$(Y_a(t))_{t\geq 0}$ obtained through the Riccati transformation and defined by 
\begin{align}\label{diffusionY}
\begin{cases} d Y_a(t) = (a - Y_a(t)^2) \, dt + dB(t) & {\rm for} \quad t\geq 0\,, \\ Y_a(0) = y \,.  \end{cases}
\end{align}
where $a \in \R$ is a fixed parameter and $y\in \R\cup\{+\infty\}$ and the diffusion $Y_a$  immediately restarts from $+\infty$ whenever an explosion occurs. 
The characterization of the marginal distribution of the minimal eigenvalue with index $k$ now reads 
\begin{equation}\label{charac.statio}
\P\left[A_k^\beta < a \right] = \P\left[ Y_a \,  \textrm{blows up to $- \infty$ at least $k+1$ times before time $L$} \right] \,.
\end{equation} 
The forthcoming study of the exit time distribution of the diffusion $Y_a$ in the small noise limit will permit us to 
analyze the (marginal) distributions of the minimal eigenvalues of the stochastic operator $\mathcal{G}_L$
in the limit $L\to \infty$ (which replaces the limit $\beta \to 0$ in the Airy case). 

In view of subsection \ref{macro.density}, we also provide a discussion in subsection \ref{transition.pdf} 
on the transition probability density of the diffusion $Y_a$, and relate this transition pdf with the limiting  
density of state of the operator $\mathcal{G}_L$.
This section corresponds to the stationary counterpart of the main study of this paper.
The problem is of course easier to solve for the Hill's operator $\mathcal{G}_L$ than for the Airy operator 
$\mathcal{H}_\beta$ thanks to the stationarity.

We will actually see later
that the law of the minimal eigenvalues of the operators $\mathcal{G}_L$ and $\mathcal{H}_\beta$ are 
similar in the respective limits  $L\to \infty$  and $\beta\to 0$.

\subsection{Definition}
In this section, 
$\P$ denotes the law of the diffusion $Y_a$ when $y = +\infty$ and $\P_y[\cdot]:=\P[\cdot|Y_a(0)=y]$ is the law
of the diffusion $Y_a$ conditionally on $Y_a(0)=y$. In particular, $\P=\P_{+\infty}$.  

We are interested in the distribution of the exit time (blowup time) $\zeta:=\inf\{t\geq 0: Y_a(t)=-\infty\}$ 
and in particular in its limit in law when $a\to +\infty$.  

\begin{figure}[h!btp] 
     \center
     \includegraphics[scale=0.7]{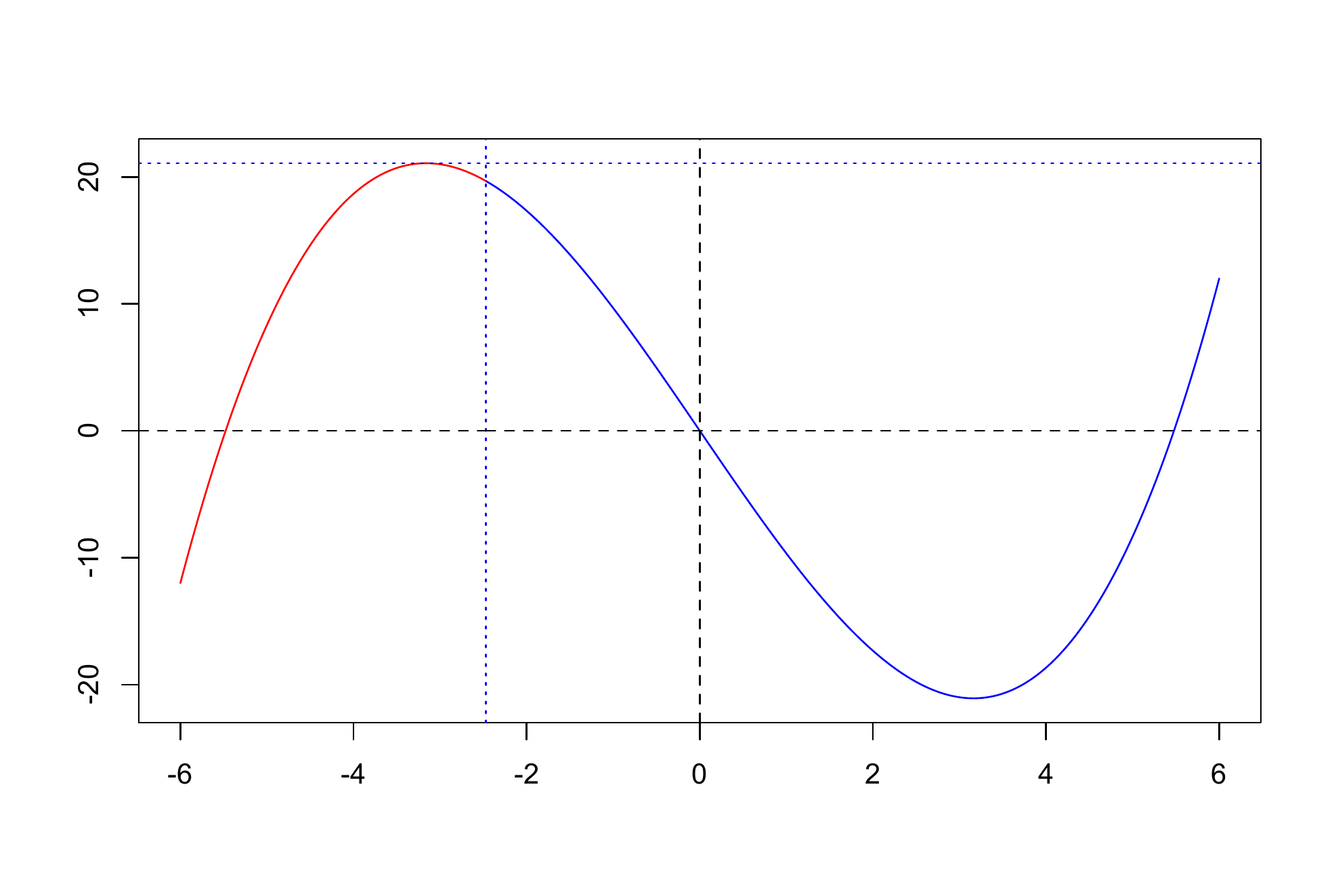}
     \caption{The potential $V(y)$ as a function of $y$. }\label{fig.potential}
\end{figure}

The diffusion evolves in a potential  $V(y) := - a y +\frac{y^3}{3}$ which
presents a local minimum in $y=\sqrt{a}$ and a local maximum in $y=-\sqrt{a}$. 
The potential barrier $\Delta V = \frac{4}{3}\, a^{3/2}$ gets very large when $a\to +\infty$ while the noise 
remains constant (see Figure \ref{fig.potential}).

If the particle starts from above the potential well (from $y=+\infty$, say) at time $0$, then 
it has to cross a very large barrier of size $\Delta V= 4\,a^{3/2}/3$ with a deep well in $\sqrt{a}$. 
From Kramer's theory \cite{kramers}, 
we expect the exit time to be distributed according to an exponential law with parameter $\sim \exp(-2\Delta V)$.  
We give a simple proof of this result by a Laplace transform method (see also \cite{texier}). 
 
\medskip
   
\subsection{Exit time distribution}

Let us introduce the Laplace transform of the first exit time $\zeta$ of the diffusion $Y_a(t)$ 
with initial position $y$ 
\begin{equation}\label{laplace.statio}
g_\alpha(y) := \E_y[e^{-\alpha \zeta}]\,. 
\end{equation}

The following proposition characterizes the Laplace transform $g_\alpha$ 
as the unique solution to a boundary value problem. The proof of this Proposition cand be found in Appendix \ref{proof.theo}.

\begin{proposition}\label{unicity.statio} 
Let $\alpha > 0$. 
Then the function $g_\alpha$ defined in \eqref{laplace.statio} is the unique bounded and twice continuously differentiable 
solution of the boundary value problem 
\begin{align}\label{ode.g}
\frac{1}{2} g_\alpha'' - (y^2-a)\, g_\alpha' = \alpha\, g_\alpha 
\end{align}
satisfying the additional boundary condition 
\begin{equation}\label{boundary.g} 
g_\alpha(y)\to 1 \quad {\rm when } \quad y \to -\infty \,.
\end{equation}  
In addition, it satisfies the fixed point equation 
\begin{equation}\label{pt.fixe.g}
g_\alpha(y) = 1 - 2\alpha \int_{-\infty}^{y} dx \int_x^{+\infty} ds \exp\left(2 \,a (s-x) + \frac{2}{3} (x^3-s^3)\right)\, g_\alpha(u)\,. 
\end{equation}
\end{proposition}

This proposition permits to derive a lot of information about the behavior of the diffusion. 
A first natural question one could ask is: what is the probability that the diffusion $(Y_a)$ explodes in a finite time? 
The answer is straightforward from Proposition \ref{unicity.statio} by taking $\alpha=0$ and noting that 
\begin{align*}
g_0(y):=\lim_{\alpha \downarrow 0}\E_{y}[e^{-\alpha \zeta}] = \P_y[\zeta< +\infty] = 1\,. 
\end{align*}  

The mean exit time of the diffusion starting at time $t=0$ from position $y$ can also be computed
explicitly from Proposition \ref{unicity.statio}. It suffices to differentiate $g_\alpha(y)$ with respect to $\alpha$ and 
then take $\alpha=0$. Eq. \eqref{ode.g} then transforms into a second order differential equation which can be solved 
explicitly. We find the mean exit time $m(a,y)$ starting from position $y$
\begin{equation}\label{eq.ma.explicit}
m(a,y) = 2 \int_{-\infty}^{y} dx \int_x^{+\infty} du \exp\left(2 a (u-x) + \frac{2}{3} (x^3-u^3)\right)\,.
\end{equation}
For $y=+\infty$, this expression $m(a,+\infty)$ (simply denoted as $m(a)$ in the sequel) 
simplifies, after two changes of variables and a further Gaussian integration, in a single integral expression 
\begin{align}\label{def.ma}
m(a)=  \sqrt{2\pi} \int_{0}^{+\infty} \frac{dv}{\sqrt{v}} \, \exp\left(2 a v - \frac{1}{6}\, v^3 \right)\,. 
\end{align} 
This explicit integral form for $m(a)$ is convenient to determine its asymptotic when $a\to +\infty$
using the saddle point method. This is done in Appendix \ref{moments.exit.time} and the estimate \eqref{asymp.ma} 
will be useful throughout the paper. 

All the moments $\E_y[\zeta^n]$ for $n\in\N$ of 
the exit time $\zeta$ (when the diffusion starts from $y$ at time $0$) can also be derived by iterating this argument (see Appendix \ref{proof.theo}). 
In particular all the moments are finite for all starting point $y$ and fixed $a$. They actually satisfy    
$\E_y[\zeta^n] \leq n!\, m(a)^n$.

\begin{rem}\label{rem.LT}
 As mentioned in the previous paragraph, it is interesting to note that we can derive a similar theorem for the Laplace transform of the first exit time of the non-stationary diffusion under consideration in the latter section $X_{\lambda}$ \eqref{eq.lang.xt}. 
Indeed, let us denote by $\zeta:= \zeta(x,t)$ the first explosion time of the diffusion process $(X_\lambda({s}))_{s\geq t}$ conditioned to start at time $t$ in $x$ and let $f_\alpha(x,t)$ for $\alpha > 0$ be its Laplace transform. For $\alpha=0$, the Laplace transform is extended by continuity 
\begin{equation}\label{laplace.in.zero}
f_0(x,t):= \lim_{\alpha \downarrow 0} f_\alpha(x,t) = \P[\zeta(x,t) < + \infty]\,. 
\end{equation} 

Similarly to Proposition \ref{unicity.statio}, we can show that the function $f_\alpha(x,t)$ is the unique solution of the following boundary value problem:
\begin{equation}\label{pde.f}
\frac{\partial f_\alpha}{\partial t} + \left( t - \lambda - x^2  \right) \frac{\partial f_\alpha}{\partial x} 
+ \frac{2}{\beta} \frac{\partial^2 f_\alpha}{\partial x^2} = \alpha \, f_\alpha\,,
\end{equation}
satisfying the additional boundary conditions 
\begin{align}\label{boundary.v1}
& f_\alpha(x,t) \rightarrow 1\quad  {\rm when} \quad x \to -\infty\quad {\rm with} \quad  t  \quad {\rm fixed} \,, \\
& f_\alpha(x,t) \rightarrow 0 \quad {\rm when} \quad x,t \to +\infty \quad {\rm together}. \label{boundary.v2}
\end{align}
 Note that this characterization is also true for $\alpha=0$: We recover the result of 
\cite[Theorem 1.7 (ii)]{bloemendal} which permits to find the cumulative distribution of the minimal eigenvalue thanks to Eq. \eqref{charac.tw} which rewrites as  
\begin{equation*}
\P[\Lambda_0^\beta > \lambda] = \lim_{x\to +\infty} \P[\zeta(x,0) = + \infty]\,. 
\end{equation*}
The generalization given here has the advantage to contain the full information about the law of the first exit time. Nevertheless, it did not permit us to derive the convergence of the point process of the explosion times.
\end{rem}

We can now establish the weak convergence when $a\to +\infty$ 
of the explosion waiting time $\zeta$ when rescaled by its mean $m(a)$ to the exponential distribution with parameter $1$. 
The result is valid for a large range of starting points, mainly all the points which are {\it in} the potential well ({\it i.e. above} the local 
maximum of the potential). 

\begin{theorem}\label{conv.exp}
Let $f:\R\to \R$ such that $a^{1/4}(f(a)+a^{1/2}) \to_{a\to +\infty} +\infty$. 
Then,
\begin{align*}
 \sup_{y\geq f(a)} \left|g_{\alpha/m(a)}(y) - \frac{1}{1+\alpha}\right| \longrightarrow 0\,.  
\end{align*}
In particular, for any $y\geq f(a)$
the first blowup time to $-\infty$ of the diffusion $(Y_a(t))$ starting from position $y$ at time $0$, rescaled as $\zeta/m(a)$, converges weakly 
when $a\to +\infty$ to an exponential law with parameter $1$. 
\end{theorem}

\begin{rem}
We recover here the prediction of Kramer's theory. Indeed
 the exit time (starting from a point inside the well) is distributed in the limit $a\to +\infty$ 
 according to an exponential law with parameter $m(a)^{-1}$, for which we have found a logarithmic equivalent 
  $m(a)^{-1}\asymp \exp(-8/3 \,a^{3/2})$ in Appendix \ref{moments.exit.time}. 
 \end{rem}

From Theorem \ref{conv.exp} and the characterization \eqref{charac.statio} of the law minimal eigenvalues of the operator 
$\mathcal{G}_L$, we get back the result due to McKean in \cite{mckean} about the fluctuations of the ground state of the Hill's operator. 
\begin{corollary}\label{conv-ground-Hill}
In the limit $L\to \infty$, the fluctuations of 
the minimal eigenvalue $A_0^L$ of the stochastic Hill operator $\mathcal{G}_L$ are governed by the Gumbel 
distribution. More precisely, we have the convergence in law as $L\to \infty$,
\begin{align*}
- 2 \cdot 3^{1/3} (\ln L)^{1/3} \left[A_0^L + \left(\frac{3}{8}\ln \frac{L}{\pi}\right)^{2/3}  \right] 
\Rightarrow e^{-x} \exp(-e^{-x}) \, dx\,. 
\end{align*}
\end{corollary}
{\it Proof.}

Using the asymptotic estimate \eqref{asymp.ma} for $m(a)$ as $a\to +\infty$, the convergence in law of the exit time 
towards the exponential distribution provided by Theorem \ref{conv.exp} 
and the characterization \eqref{charac.statio}, it is easy to see that 
\begin{align*}
\P\left[A_0^L\leq - \left(\frac{3}{8}\ln \frac{L}{\pi}\right)^{2/3}  + \frac{1}{2\cdot 3^{1/3}}  \frac{x}{(\ln L)^{1/3}} \right] \to_{L\to +\infty} 1-e^{-e^x}\,, 
\end{align*}
which yields the result. 
\qed
\subsection{Exit times point process}\label{explosions.process.statio}

Although the diffusion $(Y_a(t))$ blows up to $-\infty$ in a finite time almost surely, we can again define 
the trajectory for all time $t\geq 0$ by 
restarting the diffusion in $+\infty$ immediately after any blow up to $-\infty$. 
 
Endowed with this new definition, we introduce the empirical measure $\mu_a$ of the explosion times 
$\zeta_1<\zeta_2<\zeta_3\dots$ (with a further rescaling) defined for any Borel set $B$ of $\R_+$ as 
\begin{align}\label{mu_a}
\mu_a(B)= \sum_{i=1}^{+\infty} \delta_{\zeta_i/m(a)}(B)\,. 
\end{align}
In particular, the random variable $\mu_a([0,t])$ is the number of explosions of the diffusion in the interval $[0,m(a) t]$. 

We have seen (see Theorem \ref{conv.exp}) that the sequence of the waiting time $(\zeta_i-\zeta_{i-1})/m(a)$ 
until the next increment for $\mu_a$ converges in law  
when $a\to +\infty$ to an exponential distribution with parameter $1$.
It is therefore easy to deduce the convergence of the point process $\mu_a$ to a Poisson point process on 
$\R_+$. More precisely, we have the following result. 

\begin{theorem}\label{poisson.statio}
Let $Y_a$ the diffusion \eqref{diffusionY} with initial position $Y_a(0)=+\infty$.
The associated explosion times 
point process $\mu_a$ \eqref{mu_a} converges weakly (in the space of Radon measures on $\R_+$ equipped with 
the topology of vague convergence \cite{kallenberg}) when $a\to +\infty$ to a Poisson point process 
with intensity $1$ on $\R_+$. 
\end{theorem}

\noindent The end of this section is devoted the proof of Theorem \ref{poisson.statio}. 
\vspace{0.5cm}

\noindent {\it Proof of Theorem \ref{poisson.statio}.}

We will use the useful criterion from Kallenberg \cite{kallenberg}, which states that it is sufficient to prove that, for any finite 
union $I$ of disjoint and bounded intervals, we have 
the following convergences 
when $a\to +\infty$,  
\begin{align}
\E[\mu_a(I)] &\longrightarrow |I|\,,  \label{cv1} \\
\P[\mu_a(I) = 0] &\longrightarrow \exp(- |I| ) \label{cv2} \,,
\end{align}
where $|I|$ denotes the length of the set $I$.

Towards \eqref{cv1}, by linearity we just need to prove that $\E[\mu_a[0,t]] \to t$. 
The advantage here is that the starting point of the diffusion $Y_a$ at time $0$ is $+\infty$ so that the waiting time $\zeta/m(a)$
until the first explosion 
converges weakly to an exponential distribution according to Theorem \ref{conv.exp}. 
We have 
\begin{align*}
\E[\mu_a[0,t]] = \sum_{k=0}^{+\infty} \P\left[\mu_a[0,t] \geq k\right] =  1+ \sum_{k=1}^{+\infty} \P\left[\frac{\zeta_k}{m(a)} \leq t\right] \,, 
\end{align*}
where $\zeta_k$ is the $k$ th exit time of the diffusion $Y_a$ started at $+\infty$. 
The strong Markov property implies that, for each $k$, the random variable $\zeta_{k+1}-\zeta_{k}$ 
is independent of $(\zeta_1,\zeta_2-\zeta_1,\dots,\zeta_{k}-\zeta_{k-1})$ and has the same distribution as $\zeta_1$.
It is then easy to deduce that for any fixed $k$, $\zeta_k/m(a)$ converges in law when $a\to +\infty$ to 
the Gamma distribution $\Gamma(k,1)$ with shape and scale parameter $k$ and $1$ [Recall $\Gamma(k,1)$ is simply the law of 
a sum of $k$ independent exponential random variables with parameter $1$]. 
One can easily prove existence of a constant $C>0$ independent of $k$ such that   
$\P[\zeta_k /m(a) \leq t]\leq C/k^2$ using for instance Chebyshev's inequality since we have
$\E[(\zeta_i -\zeta_{i-1})^2/m(a)^2]\leq 1$ for all $i \leq k$ or the Cramer's large deviation principle which would give a much 
stronger bound.  
The bounded convergence theorem finally applies and gives, when $a\to +\infty$, 
\begin{align}\label{cv.expectation.mu_a}
\E[\mu_a[0,t]] \longrightarrow 1+\sum_{k=1}^{+\infty} \P\left[\Gamma(k,1) \leq t\right] = t.
\end{align}    

The second equality \eqref{cv2} is proved with the same idea. For any $k\geq 1$, using the strong Markov property 
for the diffusion $Y_a$, we easily prove that the first $k$ explosions times converge jointly in law 
to the first $k$ occurrence times $\xi_1<\xi_2<\cdots< \xi_k$ of a Poisson point process with intensity $1$ (the increments are independent and each of them converges to an exponential distribution thanks to Theorem \ref{conv.exp}) {\it i.e.}
\begin{align}\label{cv.in.law}
\frac{1}{m(a)} (\zeta_1,\zeta_2,\dots,\zeta_k) \Rightarrow_{a\to +\infty} 
(\xi_1, \xi_2,\dots,\xi_k)\,. 
\end{align}
The convergence \eqref{cv2} follows by proving first that we can consider only the first $k$ explosion times
for $k$ large enough using a large deviation argument as above (in order to bound the probability of having more than $k$ explosion 
times before time $m(a)t$) and then by using the convergence \eqref{cv.in.law} on the overwhelming event $\{\zeta_{k}>m(a)t\}$ . 
\qed

\subsection{Transition probability density}\label{transition.pdf}

We aim at describing also the probability density $p(y,t)$ of the diffusion $Y_a(t)$ at time $t$ and other related quantities. 
The main tool is the Fokker Planck equation which gives the evolution of  $p(y,t)$
and writes as 
\begin{align}\label{fokker.planck.statio}
\frac{\partial p}{\partial t} = \frac{\partial}{\partial y} \left[(y^2-a)\, p(y,t) + \frac{1}{2}\frac{\partial^2}{\partial y^2} p(y,t) \right]\,. 
\end{align}
This equation takes the form of a {\it continuity equation}. Indeed, introducing the flux 
\begin{equation*}
j(y,t) :=  (y^2-a)\, p(y,t) + \frac{1}{2}\frac{\partial^2}{\partial y^2} p(y,t)\,, 
\end{equation*}
Eq. \eqref{fokker.planck.statio} rewrites as 
\begin{equation*}
\frac{\partial p}{\partial t} = \frac{\partial j}{\partial y} \,. 
\end{equation*}
The equilibrium of the system can be characterized by finding the stationary 
solution $p_0(y)$ to the Fokker Planck equation \eqref{fokker.planck.statio} which satisfies 
\begin{equation}\label{eq.p0}
\frac{1}{2} \, p_0'(y) + (y^2-a)\, p_0(y) = J_0\,,
\end{equation}
where $J_0$ is a constant, which does not depend on $y$ or $t$. 
We can solve the ode Eq. \eqref{eq.p0} explicitly and find the constant $J_0$ using the  
additional normalization constraint $\int_\R p_0 =1$. We obtain 
\begin{equation}
 p_0(y) = 2 \, J_0(a) \,  \int_{-\infty}^y du \, \exp\left[2 a (y-u) + \frac{2}{3}(u^3-y^3)\right] \,,
\end{equation}
where  
\begin{equation}\label{flux.statio.J0}
J_0(a)= \frac{1}{m(a)} = \frac{1}{\sqrt{2\pi}} \left[ \int_{0}^{+\infty} \frac{dv}{\sqrt{v}} \, \exp\left(2 a v - \frac{1}{6}\, v^3 \right)\right]^{-1} \,. 
\end{equation}
Here we stress that $J_0(a)$ is precisely equal to the inverse of the expected exit time of
the diffusion starting from $+\infty$ \cite{lloyd,halperin}. 
%
In the limit of large $L$, 
the number $n_L$ of blow-ups to $-\infty$ of the 
diffusion $(Y_a(t))_{t\geq 0}$ during the interval $[0;L]$ (or equivalently -- using the characterization \eqref{charac.statio} --  
the number of $\mathcal{G}_L$-eigenvalues 
lying below the level $a$)
is proportional (at leading order) to the inverse of the expected exit time of the diffusion, 
\begin{equation}\label{emp.eigenvalue.statio}
n_L = L\, J_0(a) +O(\sqrt{L})\,. 
\end{equation}
We recover the formula of McKean \cite{mckean} for the limiting integrated density of states of the operator $\mathcal{G}_L$ introduced above, as 
$L\to\infty$. The scaling order $\sqrt{L}$ for the fluctuations comes from the central limit theorem.

Let us finally mention that this study corresponds to the $N=1$ case of a more general model of $N$ interacting particles. In this context, we expect the number of explosions per unit of time to display a different fluctuation order when $N$ goes to infinity because of the non-trivial interaction. See \cite{cubic} for more details about the model and this conjecture.

\section{Statistics of ${\rm SAE}_\beta$ when $\beta\to 0$}\label{sae.results}
We first state our results on the limiting marginal distributions of the minimal eigenvalues of the stochastic 
Airy operator $\mathcal{H}_\beta$  
when $\beta\to 0$. 
Then, we provide a more global information on the spectrum by computing the 
empirical eigenvalue density on the macroscopic scale.

Let us introduce a more convenient stochastic linear operator defined for $t>0$ as 
\begin{equation*}
\mathcal{L}_\beta := -\frac{d^2}{dt^2} + \frac{\beta}{4} \, t + B'(t)\,.
\end{equation*}
Denoting by $L_0^\beta<L_1^\beta< \dots< L_k^\beta\cdots$ its eigenvalues, it is easy to show that in law 
\begin{equation}\label{equality.law}
\{L_k^\beta, \, k\in \N\} \stackrel{(d)}{=} \left(\frac{\beta}{4}\right)^{2/3} \, \{\Lambda_k^\beta, \, k\in \N\}\,,
\end{equation}  
where $\Lambda_0^\beta< \Lambda_1^\beta< \dots < \Lambda_k^\beta\cdots$ are the eigenvalues of the operator $\mathcal{H}_\beta$. 
Indeed, \eqref{equality.law} follows after the change of function $\psi(t)=c\phi(t/c)$ with $c=(\frac{4}{\beta})^{1/3}$ in 
the differential equation \eqref{eq.diff.phi} satisfied by the eigenfunction $\phi$ of $\mathcal{H}_\beta$. 

The Riccati diffusion $Z_{\ell}$ associated to the stochastic linear operator $\mathcal{L}_\beta$ satisfies 
\begin{equation}\label{eq.zt}
Z_{\ell}(0) = +\infty \quad {\rm and}  \quad dZ_{\ell}(t) = \left(\frac{\beta}{4}\, t - \ell - Z_{\ell}(t)^2 \right) \, dt + dB(t)  \quad {\rm for} \quad t \geq 0
\end{equation} 
where $B$ is a Brownian motion. The law of the diffusion $Z_\ell$ is denoted $\P$ and $\P_{z,t}[\cdot]:=\P[\cdot|Z_\ell(t)=z]$ 
is the law of the diffusion $(Z_\ell(s))_{s\geq t}$ conditionally on $Z_\ell(t)=z$.

Recall that when blowing up to $-\infty$ at some time $t$, the diffusion $Z_{\ell}$ 
immediately restarts in $+\infty$ at this time $t$ and that we have the key relation 
\begin{align}\label{key.relationZ}
\P[L_k^\beta\leq \ell] =  \P\left[ Z_\ell  \mbox{ blows up to $- \infty$ at least $k$ times in $\R_+$ }\right] \,.
\end{align}

\subsection{Minimal eigenvalues of  $\mathcal{L}_\beta$} 
 
We investigate the minimal eigenvalues statistics of the linear stochastic operator $\mathcal{L}_\beta$ when $\beta \to 0$ and 
in particular the convergence of the marginal distributions of the bottom eigenvalues. 

When $\beta\to 0$, we can check that for any \emph{fixed} $\ell$ and $k$, $\P[L_k^\beta\leq \ell ]\to 1$. 
Indeed, if $\beta$ tends to $0$ while $\ell$ is fixed, the diffusion $Z_\ell$ defined in \eqref{eq.zt}
converges in law to the diffusion  
$Y_a$ of the previous section for $a=-\ell$. The probability for the diffusion $Z_\ell$ to explode $k$ times
in $\R_+$ will therefore tend to the corresponding probability for the diffusion $Y_a$, which is exactly equal to $1$. 
Recalling \eqref{key.relationZ}, 
we have the claim. 

Therefore if we look for a non trivial limit in law for the eigenvalue $L_k^\beta$ when $\beta\to 0$, the parameter $\ell$ should decrease to $-\infty$ as $\beta \to 0$ and we need to determine the rescaling of $\ell$ as a function of $\beta$. We can actually make the guess $\ell_\beta\sim -\ln(1/\beta)^{2/3}$ 
using the right tail asymptotic 
of the Tracy-Widom $\beta$ distribution \eqref{right.tail.tw} and the 
relation \eqref{equality.law} between the laws of $L_0^\beta$ and $\Lambda_0^\beta = - TW(\beta)$. 
 
In agreement with this heuristic derivation, we fix $x\in \R$ 
and set
\begin{align}\label{def.ell.beta}
\ell_\beta=\ell_\beta(x) := - \left(\frac{3}{8} \ln\frac{1}{\beta\pi}\right)^{2/3} + \frac{1}{2}\frac{1}{3^{1/3}} \, \left( \ln \frac{1}{\beta} \right)^{-1/3} x \,. 
\end{align}
We shall prove that the function $x\to \P[L_k^\beta\leq \ell_\beta(x)]$ 
converges to a non trivial cumulative distribution function on $\R$ when $\beta \to 0$. 

Similarly to the previous section, we first consider
the empirical measure $\nu_\beta$ of the explosions times $(\zeta_k)_{k\in \N}$ of the diffusion $Z_{\ell_\beta(x)}$ 
after a further (well chosen) rescaling of time. For a Borel set $A\subset \R_+$, 
\begin{align}\label{def.nu_beta}
\nu_\beta(A)= \sum_{k=1}^{+\infty} \delta_{\beta (\frac{3}{8}\ln(1/\beta))^{1/3} \, \zeta_k} (A)\,. 
\end{align}
Finally, recall the key point: for almost all $x$, the number $\nu_\beta(\R_+)$ of explosions in $\R_+$ is equal 
to the number of eigenvalues smaller than $\ell_\beta(x)$.

\begin{theorem}\label{th.conv.explosion.times}
The explosion times point process $\nu_\beta$ associated to the diffusion $Z_{\ell_\beta(x)}$  
converges weakly (in the space of Radon measures equipped with the topology of vague convergence \cite{kallenberg})
when $\beta\to 0$ to a Poisson point process with inhomogeneous intensity $e^{x} \times e^{-t} \, dt$. 

It readily implies the following convergence: for any $t < t'$, $k\in \N$,
\begin{equation*}
\P\left[\nu_\beta\left[t;  t' \right] = k \right] \longrightarrow_{\beta \to 0} \exp\left( - e^{x} \, \int_{t} ^{t'}  e^{-s}\, ds\right) \frac{\left(e^{x} \int_{t} ^{t'} e^{-s}\, ds\right)^k}{k!}\,.  
\end{equation*}
\end{theorem}

It might be useful to compare the time scale of the measure $\nu_\beta$ 
with the one used in the definition of $\mu_a$ \eqref{mu_a}.
Using the asymptotic formula Eq. \ref{asymp.ma} given in Appendix \ref{moments.exit.time}, we can check that
\begin{equation*}
m(-\ell_\beta(0))^{-1} = m\Big( \big(\frac{3}{8} \ln(1/(\beta\pi))\big)^{2/3}\Big) ^{-1} \sim_{\beta\to 0} \beta \left(\frac{3}{8}\ln \frac{1}{\beta}\right)^{1/3}.
\end{equation*}  
In such a way, the two time scales are actually the same up to exchanging the parameters $a$ and 
$-\ell_\beta(0)=(\frac{3}{8} \ln\frac{1}{\beta\pi})^{2/3}$.
More generally, note that for any $x\in \R, t \in \R_+$,
\begin{align*}
m\left(-\ell_\beta(x-t))\right) ^{-1} \sim_{\beta \to 0} \beta \left(\frac{3}{8}\ln \frac{1}{\beta}\right)^{1/3} e^{x-t}\,.
\end{align*}
The scenario described in Theorem \ref{th.conv.explosion.times} is therefore fairly simple. 
When $\beta\to 0$, the diffusion $Z_{\ell_{\beta}(x)}$ feels the evolution in time due to the linear term $\frac{\beta}{4}t$ in the drift
but in a rather trivial way:  
setting 
\begin{align*}
s :=\frac{t}{\beta (\frac{3}{8} \ln(1/\beta))^{1/3}} 
\end{align*}
to work in the appropriate time scale, 
the explosion times process of $Z_{\ell_\beta(x)}$ is somehow the same as the one of a stationary ``frozen" system  
which evolves with time $t$ {\it adiabatically} 
such that the parameter $a$ in the drift of the diffusion $Y_a$ evolves slowly with time $t$ as
\begin{align}\label{def.at}
a:= - \ell_\beta(x) +\frac{\beta}{4}s = -\ell_\beta(x-t)\,. 
\end{align}

\medskip

Finally let us come to the initial motivation of this paper which comes from interesting
questions, recently asked in the literature in \cite{jp-alice,satya,johansson}, about the existences and characterizations 
of possible distributions in extreme value theory which would 
interpolate between the Tracy-Widom $\beta=1,2,4$  
(maximum of {\it highly} correlated random variables) and the Gumbel law (maximum of {\it weakly} correlated random variables). 

The following Theorem, obtained as a straightforward application of Theorem \ref{th.conv.explosion.times},  
establishes and describes precisely the progressive deformation of the Tracy-Widom $\beta$ laws into a Gumbel 
law when $\beta \to 0$.  

\begin{theorem}\label{TracyWidom.Gumbel}
When properly rescaled and centered, the Tracy-Widom $\beta$ law converges weakly to the Gumbel law. 
More precisely, when $\beta\to 0$, the following convergence in law holds
\begin{equation*}
2 \cdot 3^{1/3} \cdot \left( \ln \frac{1}{\beta} \right)^{1/3} \left[ \left(\frac{\beta}{4}\right)^{2/3} \, TW(\beta)-  \left(\frac{3}{8}\right)^{2/3}\, \left(\ln\frac{1}{\beta\pi}\right)^{2/3} \right]\Longrightarrow e^{-x} \exp(-e^{-x}) \, dx\,.  
\end{equation*}
\end{theorem}

Comparing Theorem \ref{TracyWidom.Gumbel} to Corollary \ref{conv-ground-Hill}, we see that the fluctuations of the ground state 
of the stochastic Hill and Airy operators are very similar, governed by the Gumbel 
distribution, with the same scalings under the relation $L=1/\beta$.

The convergences of the marginal distributions of the other minimal eigenvalues with general index $k$ can also be deduced 
from Theorem \ref{th.conv.explosion.times} which implies that
\begin{align}\label{other.vps.cv}
\P\left[L_k^\beta\leq \ell_\beta(x)\right] = \P\left[\nu_\beta[0;+\infty] \geq k+1\right]  \rightarrow_{\beta\to 0}
1- e^{-e^{x}} \, \sum_{i=0}^k \frac{(e^{x})^i}{i!}\,. 
\end{align}
We can check that the function of $x$ in the right hand side is indeed a cumulative distribution function 
and the weak convergence of $L_k^\beta$ follows.
Note also that the right hand side of \eqref{other.vps.cv} 
corresponds exactly to the cumulative distribution of the point with index $k$ 
of a Poisson point process on $\R$ with intensity $e^{x} dx$. 

It is therefore natural to conjecture that the empirical measure $\rho _\beta$ of the eigenvalues of $\mathcal{L}_\beta$
defined, for any $A$ Borel set of $\R$, as 
\begin{align}\label{density.vp}
\rho_\beta(A) = \sum_{k=0}^{+\infty} \delta_{L_k^\beta}(A)\,,
\end{align} 
converges weakly when considered on the {\it microscopic scaling region of the
minimal energies} as in this paragraph to a Poisson point process on $\R$ 
with intensity $e^{x}dx$, when $\beta \to 0$.

Let us finally consider the matching between the microscopic (where one zooms in the measure $\rho_\beta$ 
on the microscopic region of the bottom eigenvalues) and the macroscopic regime where the sets $A\subset \R$ remains fixed as $\beta \to 0$. 

In section \ref{macro.density}, we compute asymptotically when $\beta\to 0$ the empirical spectral measure on the macroscopic scale. 
Our result reads 
\begin{align}\label{asymp.rho.beta}
\rho_\beta(A) = \frac{4}{\beta} \int_A J_0(\ell) d\ell\, + O(1)
\end{align}
where $O(1)$ is a constant of order $1$ as $\beta\to 0$. In particular, $\rho_\beta(]-\infty;\ell])$ is the number of eigenvalues below the level $\ell$. 
Using the result obtained in Appendix \ref{consistency.pt}, --- where we 
check that the second order term $O(1)$ in \eqref{asymp.rho.beta} remains negligible compared to the leading order when $\beta \to 0$ and
$A=]-\infty; \ell_\beta(x)]$--- we see that for any fixed $x$, 
\begin{align*}
\rho_\beta\left(]-\infty;\ell_\beta(x)]\right) \sim_{\beta\to 0} e^{x}\,,  
\end{align*}
where we have also used the asymptotic \eqref{intJ0.dl}. 

Therefore, if one considers a Poisson point process $\mathcal{P}$ 
on $\R$ with intensity $4 J_0(\cdot)/\beta$, then the probability that 
the $k+1$ bottom points of $\mathcal{P}$ 
are below the level $\ell_\beta(x)$ converges when $\beta\to 0$ as
\begin{align*}
\P\left[\mathcal{P}\left(]-\infty; \ell_\beta(x)]\right) \geq k+1 \right] \longrightarrow 1- e^{-e^{x}} \, \sum_{i=0}^k \frac{(e^{x})^i}{i!}\,. 
\end{align*}
The matching between the two microscopic and macroscopic regimes, 
if the Poisson point process convergence be true, would be smooth.

Let us mention two articles \cite{molcanov1,molcanov2} from the literature on random Schrodinger operators (with \emph{stationary} potential)
which establish that the bottom eigenvalues have Poissonian statistics in the limit of infinite support  
for different kind of potential, including sums of Dirac masses and smooth 
functions of diffusions. 

\subsection{Proof of Theorem \ref{th.conv.explosion.times}.} \label{main.proof}

Let us fix $x \in \R$ and denote simply by $Z$ the diffusion $Z_{\ell_\beta(x)}$ in this proof. From Kallenberg's theorem \cite{kallenberg}, we just need to see that, for any finite 
union $I$ of disjoint and bounded intervals, we have when $\beta\to 0$,  
\begin{align}
\E[\nu_\beta(I)] &\longrightarrow e^x \int_I e^{-t}\, dt \,,  \label{cv1bis} \\
\P[\nu_\beta(I) = 0] &\longrightarrow \exp(- e^x \int_I e^{-t}\, dt ) \label{cv2bis} \,. 
\end{align}
Denote by $[t_1;t_2]$ the right most interval of $I$,  
by $J$ the union of disjoint and bounded intervals such that $I=J\cup [t_1;t_2]$ and by $t_0$ the supremum of $J$.
Note that $t_0< t_1$. 

To simplify notations, set $s_i:= t_i/(\beta (\frac{3}{8}\ln\frac{1}{\beta})^{1/3}) $ for $i=0,1,2$ (the times $s_i$ will correspond to the real time scale for the diffusion $Z$).

Thanks to the linearity of the expectation, it is enough to prove \eqref{cv1bis} for intervals $I$ of the form $I=[0,t]$. 

On the other hand for \eqref{cv2bis}, the simple Markov property yields   
\begin{align*}
\P[\nu_\beta(I) = 0]  = \P\left[1_{\nu_\beta(J) = 0} \P_{Z(s_0) ,s_0} \left[\nu_\beta[t_1;t_2]=0 \right] \right] \,. 
\end{align*}
Showing \eqref{cv2bis} therefore reduces to see that with probability going to $1$ as $\beta\to 0$, 
\begin{equation*}
\P_{Z(s_0) ,s_0} \left[\nu_\beta[t_1;t_2] =0\right] \longrightarrow _{\beta\to 0} \exp(- e^x \int_{t_1}^{t_2} e^{-t}\, dt )\,. 
\end{equation*}

For both \eqref{cv1bis} and \eqref{cv2bis}, the idea is to decompose the interval $[s_1;s_2]$ (with $s_1=0$ for \eqref{cv1bis}) 
into a finite number of small
intervals of length $\delta := \eps /( \beta (\frac{3}{8}\ln(1/\beta))^{1/3})$ 
and to approximate the number of explosions of $Z$ on each small interval of the subdivision by those of stationary diffusions, thanks to the {\it increasing property}. 

We will use in fact a random subdivision instead of a deterministic one to avoid technical issues due to special points.
To this end, let us define a sequence $(\tau_k)_{k \in \N}$ of i.i.d. random variables with uniform law in $[0;1]$, 
independent of the diffusion $Z$. 
Let $\delta$ small enough such that $0 <\delta < s_1-s_0$.
Then, we construct iteratively the sequence of random times $(S_k)_{k\geq 0}$
such that $S_0 := s_1 - \delta \tau_0$ \footnote{In the case where $s_0=0$, we simply set $S_0:=0$.}, 
$S_k := S_{k-1} + \delta \tau_{k}$ for $k \geq 1$ 
and the stopping times defined as:
\begin{align*}
K_1 &:= \inf\{k \geq 1 \;:\; S_k \geq s_1\}\,,\\
K_2 &:= \inf\{k \geq 1 \;:\;S_k  \geq s_2 \}\,.
\end{align*}
In the following ${\P}$ refers to expectation with respect to both $Z$ and the $\tau_k$.

\begin{figure}[!h]
 \centering
\includegraphics[width = 8cm]{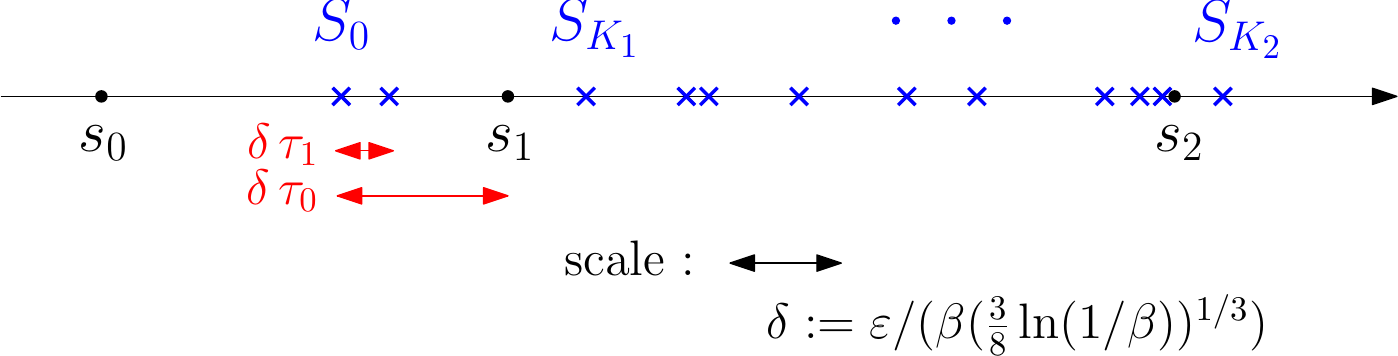}
\caption{Definition of the times $(S_k, k\geq 0)$.} 
\end{figure}

Set $a := - \ell_\beta(x)$ such that the drift of $Z$ at time $s$ simply equals $a + \beta s/4 - Z(s)^2$.
On each interval $[S_k,S_{k+1}]$ of the random subdivision, we define the following two diffusions $m_k$
and $M_k$ (independent of the times $\tau_k$) driven by the same Brownian motion as $Z$. 
\begin{align*}
&m_k (S_k) = Z(S_k)\,, \quad \mbox{and} \quad dm_k(s) = 
\left( a + \frac{\beta}{4} S_k \,  - m_k^2(s)\right) ds + dB_s \mbox{ for } s\in [S_k;S_{k+1}] \,,\\
&M_k(S_k) = Z(S_k) \,, \quad \mbox{and} \quad  dM_k(s) = \left( a + \frac{\beta}{4} S_{k+1} - M_k^2(s)\right) ds + dB_s \mbox{ for } s\in [S_k;S_{k+1}]\,. 
\end{align*}
The {\it increasing property} implies that the number of explosions 
$\nu_\beta[S_k;S_{k+1}]$ of the diffusion $Z$ is stochastically dominated from 
above (respectively from below) by the number of explosions of the diffusion $m_k$ (respectively $M_k$).  
Indeed, the drifts of the diffusion $m_k(s), Z(s)$ and $M_k(s)$ are in increasing order for $s\in [S_k;S_{k+1}]$
\begin{equation*} 
a + \frac{\beta}{4} S_k \leq a + \frac{\beta}{4} s \leq a + \frac{\beta}{4} S_{k+1}\,. 
\end{equation*}

The drifts of the diffusions $m_k$ and $M_k$ do not depend on time $s$  
so that the previous section applies on each $[S_k;S_{k+1}]$. 

Let us introduce the so-called downside region $D$ and its complementary upside one $U$ we define as 
\begin{align}\label{regionLa}
D := \{y \in \R: y \le -a^{1/2}+ \frac{(\ln a)^{1/4}}{a^{1/4}}\} \,, \quad U := \R \setminus D\,. 
\end{align}
The choice $f(a)=  (\ln a)^{1/4}$ is large enough for our purposes (it tends to infinity) and small enough such that the process will spend little time in the region $D$ (see Lemma \ref{lem.La_t}).

From the previous section, we know that \emph{as long as the starting point $Z(S_k)$ belongs to the region $U$} (which indeed is contained in the interval where the convergence applies),
the two respective explosions point processes of the diffusions $m_k$ and $M_k$ 
converge weakly in the space of Radon measure, when the time scale is renormalized 
respectively by $m(a + \beta S_k/4)$ and $m(a + \beta S_{k+1}/4)$,
to Poisson point processes with intensity $1$, independently of the exact location of $Z(S_k)$. 

Therefore
we need to prove that 
\begin{equation}\label{proba.union}
\mathcal{C}:=\bigcap_{k=0}^{K_2-1}\{Z(S_k) \in U\}
\end{equation} 
has a probability going to $1$ when $\beta \to 1$.

We first prove that the occupation time of the region $D$
by the particle $Z$ on the region of times considered tends to $0$ in expectation. 
This is the content of the following Lemma whose proof is deferred at the end of this section. 

\begin{lemma}\label{lem.La_t}
Fix $t_0'<t_1' < t_2'$ and as before $s_i':=t_i'/[\beta\, (\frac{3}{8}\ln\frac{1}{\beta})^{1/3}],i=1,2$.  
Define 
\begin{align*}
T_\beta(s_1',s_2') := \int_{s'_1}^{s'_2} 1_{\{Z(u) \in D\}} \, du \,. 
\end{align*}
Then, for any $z\in \R$ there exists $\eta>0$ and a positive constant $C$ 
independent of $\beta$ such that for all $\beta>0$,
\begin{align*}
\E_{z,s_0'}[T_\beta(t_1',t'_2)] \leq C \, \Big(\ln \frac{1}{\beta}\Big)^{-\eta} \,.
\end{align*}
\end{lemma}
We are now ready to prove that $\P(\mathcal{C}) \to 1$.
Indeed, using the notations introduced above,
the probability of $\mathcal{C}^c$ is easily bounded from above by
\begin{align}
{\P}_{Z(s_0) ,s_0} &\left[\bigcup_{k=0}^{K_2-1}\{Z(S_k) \in D\}\right] \leq 
\sum_{k=0}^{4(t_2-t_1)/\eps} {\P}_{Z(s_0) ,s_0} \left[Z(S_k) \in D \right]+ 
{\P}\big[\sum_{k=1}^{4(t_2-t_1)/\eps} \tau_k \leq  \frac{t_2-t_1}{\eps} \big]\,. \label{union} 
\end{align} 
The second probability on the RHS of \eqref{union} is easily bounded by noting simply that the $\tau_k$ have mean $1/2$ so that 
the empirical sum of the $\tau_k$ should be of order $2(t_2-t_1)/\eps$. 
From large deviation theory, we know that this probability is of order $\exp(-c/\eps)$ 
where $c$ is a positive constant independent of $\eps$. 

Thanks to the independence between the diffusion 
$Z$ and the sequence $(S_k)$, the sum in \eqref{union} can be bounded as 
\begin{align*}
\sum_{k=1}^{4(t_2-t_1)/\eps} {\P}_{Z(s_0) ,s_0}\left[Z(S_k) \in D \right]&= 
\sum_{k=1}^{4(t_2-t_1)/\eps} {\E}_{Z(s_0) ,s_0}[\int_{S_{k-1}}^{S_{k-1}+\delta} 1_{\{Z(\tau) \in D \}}\, \frac{ d\tau}{\delta} ] \\
&\leq  2\beta (\frac{3}{8}\ln\frac{1}{\beta})^{1/3}\varepsilon^{-1} \, {\E}_{Z(s_0) ,s_0}\left[\int_{S_0}^{S_{4(t_2-t_1)/\eps}}1_{\{Z(\tau) \in D\}} \, d\tau \right]\,. 
\end{align*}
Using again the Cramer's theorem of large deviation theory, we 
can show that, with probability going to $1$ when $\eps\to 0$, the empirical sum 
$S_{4(t_2-t_1)/\eps} - S_0= \delta \sum_{k=1}^{4(t_2-t_1)/\eps} \tau_k$  is smaller than $\delta+ 3(s_2-s_1)$ 
\footnote{the choice of $3$ is arbitrary.}.
Lemma \ref{lem.La_t}, applied for fixed times $s'_0:=s_0 <s'_1<s_2'$ such that
$[s_1-\delta; s_1+ 3(s_2-s_1)]\subset [s_1';s_2']$ and for $z=Z(s_0)$, permits to 
conclude that there exist $\eta >0$ and two positive constants $c$ and $C>0$ 
both independent of $\beta,\eps, \eta$ such that 
\begin{align*}
{\P}_{Z(s_0) ,s_0}&\left[\bigcup_{k=1}^{K_2-1}\{Z(S_k) \in D\}\right] \leq 
C \beta \, (\ln \frac{1}{\beta})^{1/3-\eta} +\, e^{-c/\eps}\,. 
\end{align*} 

We are now ready to prove \eqref{cv1bis} and \eqref{cv2bis}. 
To simplify notations, set $N_k$, 
$N_k^-$ and $N_k^+$ for the number of explosions 
of the diffusions $Z, M_k$ and $m_k$ in the interval $[S_k;S_{k+1}]$. Recall that a.s. for all $k$, $N_k^- \leq N_k\leq N_k^+$

We begin with \eqref{cv1bis}. Similarly to the stationary case, we could take $s_1=0$ in this case for which the situation is slightly easier. 
Indeed, at the starting point $s_1=0$, $Z(s_1)=+\infty\in U$ and we do not need to introduce an independent random time 
just before $s_1$ to start from a nice position. Nevertheless, we treat the case $s_1>0$ to have consistent notations with the ones required in the proof 
of \eqref{cv2bis}. 
 
Thanks to the increasing property, we can bound from above the mean number of explosions of $Z$ as
\begin{align}
\E[\nu_\beta[t_1,t_2]] \leq \E\left[\sum_{k=0}^{K_2-1}\E\left[ N_k^+ \left| Z(S_k)\in U\right.\right ] \right]+ 
\E[\nu_\beta[t_1,t_2] 1_{\mathcal{C}^c}]\,. \label{upper.bound.expect}
\end{align} 
The second expectation in \eqref{upper.bound.expect} is negligible.
To show this, we use the Cauchy Schwartz inequality and the increasing property which permits to bound 
the second moment of $\nu_\beta[0,t_2]$ by the second moment of the number 
of explosions in the interval $[s_1,s_2]$ of the diffusion $Y_{a+ \beta s_1/4}$ studied in section \ref{exit.time.statio}. 
It is easily seen that this second moment is bounded independently of $a$ thanks to Cramer's theorem for example as was done 
to bound the expectation of the number of explosions in the proof of Theorem \ref{poisson.statio}. 

Gathering the above arguments, we can deduce the convergence of the first term.
Indeed, conditionally on $\{Z(S_k) \in U\}$, the explosion times process  
of $m_k$ converges when time is renormalized by $m(a + \beta S_{k}/4)$, to a Poisson point process with intensity $1$. Here we must 
pay attention to the renormalization of time:  
Noting that 
\begin{align*}
a + \beta S_k/4 = \left(\frac{3}{8} \ln\frac{1}{\beta\pi}\right)^{2/3} + \frac{1}{2}\frac{1}{3^{1/3}} \, \left( \ln \frac{1}{\beta} \right)^{-1/3}
\left( t_1 -x + \eps (- \tau_0+\sum_{i=1}^{k} \tau_i ) \right)\,, 
\end{align*}
we can find the limit of the time scales ratio
\begin{align*}
\beta \left(\frac{3}{8} \ln \frac{1}{\beta \pi}\right)^{1/3} \, m(a + \beta S_k/4) \longrightarrow_{\beta \to 0} \exp(t_1 -x + \eps (-\tau_0+ \sum_{i=1}^{k} \tau_i ) )\,.  
\end{align*}
Therefore, when $\beta\to 0$, we have
\begin{align*}
\E\left[ N_k^+ \left| Z(S_k)\in U\right.\right ]  \longrightarrow 
\E[\tau_{k+1} \exp(x-t_1-\eps (-\tau_0+ \sum_{i=1}^{k} \tau_i ) )]\,, 
\end{align*}
which finally gives taking $\eps \to 0$ and thanks to the convergence of the Riemann sum associated to the subdivision $\tau_1< \tau_1+\tau_2<\cdots< \tau_1+\dots+\tau_k<\cdots$ that
\begin{align}\label{conv.expect.1}
\limsup_{\beta\to 0} \E[\nu_\beta[t_1,t_2]] \leq  e^{x} \int_{t_1}^{t_2} e^{-u}\, du \,.
\end{align}

\medskip

\noindent For the lower bound, 
\begin{align*}
&\E[\nu_\beta[t_1,t_2]]= \E\left[\sum_{k=0}^{\infty}\E\left[ N_k^- \, \, 1_{\{s_1\leq S_k \leq S_{k+1} \leq s_2\}} 1_{\{ Z(S_k)\in U \}} \right ] \right] \\
& \geq \E\left[\sum_{k=0}^{4(t_1-t_0)/\eps }\E\left[ N_k^- \, \, 1_{\{s_1\leq S_k \leq S_{k+1} \leq s_2\}} \left| Z(S_k)\in U\right.\right ] \right]\\ 
&- \E\left[\sum_{k=0}^{4(t_1-t_0)/\eps } \E\left[ N_k^- \, \, 1_{\{s_1\leq S_k \leq S_{k+1} \leq s_2\}}\left| Z(S_k)\in 
U\right.\right ] 
\P\left[Z(S_k)\in D\right]\right] \,.
\end{align*}
It suffices to show that the second term is negligible and apply the same method as was done above to obtain \eqref{conv.expect.1}.  
The second term is bounded from below by 
\begin{align}
- \E\left[\mbox{Number of explosion of $Y_{a + \beta s_1/4}$ in } [s_1,s_2]\right] \times \sup_{0\leq k \leq \frac{4(t_1-t_0)}{\eps}}  
\frac{\P[Z(S_k)\in D] }{1- \P[Z(S_k)\in D]} \,. \label{le.sup}
\end{align}
Lemma \ref{lem.La_t} applied for $t'_0 := 0<t'_1<t'_2$ such that $t'_1<t_1-\eps$ and $t'_2> t_1+ 4(t_2-t_1)$ gives a uniform (independent of $k \leq 4(t_2-t_1)/\eps$) bound of
$\P[Z(S_k)\in D]$ which goes to $0$ as $\beta\to 0$. 
It easily follows that the quantity \eqref{le.sup} tends to $0$ (we have already seen that the first term is bounded independently of $a$).

We now turn to the convergence \eqref{cv2bis}. By the increasing property, we have 
\begin{align*}
{\P}&_{Z(s_0) ,s_0}\left[\nu_\beta[t_1;t_2] =0\right] \\ 
&\leq 
{\P}_{Z(s_0) ,s_0}\left[ \prod_{k=K_1}^{K_2-2} 
\P_{Z(S_k) ,S_k} \left[ N_k^-=0 \left| Z(S_k) \in U, (\tau_k)_{k\in \N}  \right. \right]\right] + 
{\P}_{Z(s_0) ,s_0}[\mathcal{C}^c]\,.
\end{align*}

Using the previous arguments, we get
\begin{align*}
\limsup_{\beta\to 0}\, {\P}_{Z(s_0) ,s_0}\left[\nu_\beta[t_1;t_2] =0\right] \leq 
{\E}\left[ \prod_{k=K_1}^{K_2-2} \exp\left(- \eps \tau_{k+1} \exp(x-t_1-\eps (-\tau_0+ \sum_{i=1}^{k+1} \tau_i ) ) \right) \right] + e^{-c/\eps}\,. 
\end{align*}
Taking the limit $\eps\to 0$, we obtain as before
\begin{align*}
\limsup_{\beta\to 0} {\P}&_{Z(s_0) ,s_0}\left[\nu_\beta[t_1;t_2] =0\right]  \leq \exp(- e^{x} \int_{t_1}^{t_2} e^{-u} \, du )\,. 
\end{align*} 
 
For the lower bound, we use again the increasing property which allows to write 
 \begin{align*}
&{\P}_{Z(s_0) ,s_0}\left[\nu_\beta[t_1;t_2] =0\right]   \\ \geq
&{\P}_{Z(s_0) ,s_0}\left[ \prod_{k=0}^{K_2-1} 
\P_{Z(S_k) ,S_k} \left[N_k^+=0 \left| Z(S_k) \in U, (\tau_k)_{k\in \N}  \right. \right]\right] - 
{\P}_{Z(s_0) ,s_0}[\mathcal{C}^c]\,.
\end{align*}
Using the same arguments as above, we can deduce that 
\begin{align*}
\liminf_{\beta\to 0} {\P}&_{Z(s_0) ,s_0}\left[\nu_\beta[t_1;t_2] =0\right]  \geq \exp(- e^{x} \int_{t_1}^{t_2} e^{-u} \, du )\,, 
\end{align*} 
 and \eqref{cv2bis} is proved. 
 
{\it Proof of Lemma \ref{lem.La_t}.}

The key estimate to prove this Lemma is given by \cite[Proposition 10]{laure}, which is recalled in Appendix \ref{proof.paper.laure}. Let us recall it here.
We denote by $\zeta^{(u)}:=\inf\{v\geq 0: Z(u+v)=-\infty\}$ the waiting time of the first explosion after $u$.
There exists a constant $c>0$ independent of $\beta$
such that for all $u \in [s'_1;s'_2]$ \footnote{$[s'_1;s'_2]$ indeed is contained in the region where Lemma \ref{lemma.laure} applies.}
\begin{equation}\label{nice.estimate}
\P_{-a^{1/2}+ \frac{(\ln a))^{1/4}}{a^{1/4}}, u}\left[\zeta^{(u)} < \frac{\ln a }{\sqrt{a}} \right] \geq \exp(-c\sqrt{\ln a})\,. 
\end{equation}
We can rewrite the inequality \eqref{nice.estimate} in the following way, using the asymptotic expansion of $a$ as a function of $\beta$:
\begin{equation}\label{nice.estimate2}
\P_{-a^{1/2}+ \frac{(\ln a)^{1/4}}{a^{1/4}}, u}\left[\zeta^{(u)} < C\,(\ln \frac{1}{\beta})^{-\eta} \right] \geq \exp(-c\, (\ln \frac{1}{\beta})^{1/3})\,. 
\end{equation}
Eq. \eqref{nice.estimate2} gives a lower bound of the probability 
for the diffusion $Z$ starting from the right most point of the interval $D$ at time $u$ 
to explode to $-\infty$ 
in a short time. 
For the sake of completeness, we rewrite the proof of the estimate Eq. \eqref{nice.estimate} in Appendix \ref{proof.paper.laure}
with more details than originally given in \cite{laure}.  

The idea to obtain the Lemma is to use \eqref{nice.estimate2} 
to relate between the time spent in the region $D$ by the process $Z$
and the number of its explosions which is bounded on the interval $[s'_1;s'_2]$. In the following, the probability $\P$ is
with respect to the diffusion starting from $z$ at time $s'_0$.  

From \eqref{nice.estimate}, we can control for any given $u>0$, the probability that the diffusion $Z$
is in the region $D$ at time $u$: Indeed, we have  
\begin{align}
&\P[Z(u) \in D] = \P\left[Z(u) \in D, \zeta^{(u)} \leq C\,(\ln \frac{1}{\beta})^{-\eta} \right] + 
\P\left[Z(u) \in D, \zeta^{(u)} > C\,(\ln \frac{1}{\beta})^{-\eta}  \right]\notag \\
&\leq \P\left[Z(u) \in D,\zeta^{(u)} \leq C\,(\ln \frac{1}{\beta})^{-\eta}  \right] + \left(1- \exp(-c\, (\ln \frac{1}{\beta})^{1/3}) \right) \, 
\P\left[Z(u) \in D \right] \label{trick.proba1}
\end{align}
where we have used the simple Markov property and the increasing property in the second line. 
Eq. \eqref{trick.proba1} readily rewrites as 
\begin{align}
&\P\left[Z(u) \in D \right]  \leq\exp(c\, (\ln \frac{1}{\beta})^{1/3}) \, \P\left[Z(u) \in D, \zeta^{(u)} \leq C\,(\ln \frac{1}{\beta})^{-\eta} \right]\notag\\ 
& \leq\exp(c\, (\ln \frac{1}{\beta})^{1/3})\, \P\left[\mbox{The interval } \left[u,u +  C\,(\ln \frac{1}{\beta})^{-\eta} \right] \label{trick.proba2}
\mbox{ contains at least one explosion}\right]\,.
\end{align}
Note that, denoting by $k=\nu_\beta[t'_1,t'_2]$ the (random) number of explosions in the interval $I_{t,t'}$
and as before by $0< \zeta_1<\zeta_2 < \dots
<\zeta_k$ the explosion times, we have almost surely 
\begin{equation*}
\int_{s'_1}^{s'_2} {1}_{\{\exists i : \zeta_i\in[u,u+C\,(\ln \frac{1}{\beta})^{-\eta}  ] \}}  du \leq C\,(\ln \frac{1}{\beta})^{-\eta}  \, \nu_\beta[t'_1,t'_2]\,. 
\end{equation*} 
Therefore, integrating Eq. \eqref{trick.proba2} with respect to $u$ in the interval $I_{t,t'}$, we finally obtain the inequality 
\begin{align*}
\E[T_\beta(s'_1,s'_2)] \leq C\,(\ln \frac{1}{\beta})^{-\eta} \,  \exp(c\, (\ln \frac{1}{\beta})^{1/3})  \,  \,  
\E\left[\nu_\beta[t'_1,t'_2] \right]\,. 
\end{align*} 
In addition we easily check, using the increasing property and the convergence \eqref{cv.expectation.mu_a} 
proved in the previous section, that the mean number of explosions $\E[\nu_\beta[t'_1,t'_2]]$ 
is bounded independently of $\beta$. The Lemma follows.

\qed

\subsection{Empirical spectral measure on the macroscopic scale}\label{macro.density}
In this subsection, we derive the empirical spectral measure $\rho_\beta$ defined in \eqref{density.vp}
of the stochastic linear operator $\mathcal{L}_\beta$ on the macroscopic scale, i.e. {\it without zooming} in the minimal eigenvalues scaling region.

As mentioned above, we expect the number $\rho_\beta(]-\infty;\ell])$, which is almost surely finite for any $\beta>0$ (see \cite{virag}), to tend 
to $+\infty$ when $\beta\to 0$. 
The minimal eigenvalues are indeed going to $-\infty$ as $\beta\to0$ and if one sets $\beta = 0$ abruptly so that 
the linear confining term disappears, then the operator $\mathcal{L}_0:= -\frac{d^2}{dt^2}+B'_t$
has an infinite number of eigenvalues below any given level $\ell\in \R$ (see subsection \ref{exit.time.statio}).
The spectral statistics of this operator
have in fact been extensively studied 
in the literature (see \cite{lloyd,halperin,mckean,texier}). 
For this study, the eigenfunctions are restricted 
to a finite interval $[0;L]$ with Dirichlet boundary conditions in $t=0$ and $t=L$ 
so that the number of eigenvalues below a certain level remains finite. 
The authors investigate the spectrum of $\mathcal{L}_0$ in the large $L$ limit. 
The minimal eigenvalues statistics of the operator $\mathcal{L}_\beta$  
when $\beta \to 0$ as described above are different from the statistics found in \cite{mckean,texier}.
We shall see that the limiting empirical eigenvalue density differs as well  
(to be compared with the density found in \cite{lloyd,halperin}).

To compute $\rho_\beta(]-\infty;\ell])$, we use again the fact that the number of $\mathcal{L}_\beta$-eigenvalues strictly less than $\ell$ equals 
the number of explosions of the diffusion $(Z_t)_{t\geq 0}$ (defined in Eq. \eqref{eq.zt}) on $\R_+$. 
Counting the blow-ups along the trajectory $(Z_t)_{t\geq 0}$ can be done computing the 
flux of particles at $z=-\infty$ in the system, i.e. the number of particles going through $-\infty$ 
per unit of time. In this non stationary system, 
the flux depends on time $t$ and on position $z$ and will be denoted $J_\beta(z,t)$ in the following.

As the particle immediately restarts in $z=+\infty$ when blowing up to $-\infty$, the flux in $z=-\infty$ is equal 
to the flux in $z=+\infty$. 
We can again write the Fokker Planck equation which gives the evolution 
of the transition probability density $q_\beta(z,t)$ of the diffusion $Z_t$ as
\begin{equation}\label{eq.FPns}
\frac{\partial}{\partial t} q_\beta(z,t) = \frac{\partial}{\partial z} \left[(z^2+\ell-\frac{\beta}{4} t ) q_\beta(z,t)+ \frac{1}{2} \frac{\partial}{\partial z} q_\beta(z,t)\right]  \,.
\end{equation}
This equality can in turn be rewritten as a continuity equation 
\begin{equation*}
\frac{\partial q_\beta}{\partial t} = \frac{\partial J_\beta}{\partial z} 
\end{equation*}
where $J_\beta (z,t)=(z^2+\ell-\frac{\beta}{4} t ) q_\beta(z,t)+ \frac{1}{2} \frac{\partial}{\partial z} q_\beta(z,t)$ is the flux 
of probability in $z$ at time $t$.
By definition of the flux, the number 
of explosions to $-\infty$ of the diffusion $(Z_t)_{t\geq 0}$ on $\R_+$ is
\begin{equation*}
\rho_\beta(]-\infty;\ell])= \int_0^{+\infty} J_\beta(-\infty,t) \, dt\,. 
\end{equation*}
Re scaling time as before by setting $p_\beta(z,t):= q_\beta(z,\frac{4}{\beta} t)$, Eq. \eqref{eq.FPns} 
becomes 
\begin{equation}\label{eq.fp.pbeta}
\frac{\beta}{4} \frac{\partial}{\partial t} p_\beta(z,t) = \frac{\partial}{\partial z} \left[(z^2+\ell- t ) \, p_\beta(z,t)+ \frac{1}{2} \frac{\partial}{\partial z} p_\beta(z,t) \right]  \,.  
\end{equation}
We now use perturbation theory in the limit $\beta\rightarrow 0$ in Eq. \eqref{eq.FPns} in order to obtain an approximation
for $p_\beta$ for small $\beta$. 
The method consists in searching a solution of Eq. \eqref{eq.FPns} valid at small $\beta$ under the form
\begin{equation}\label{eq.pt.pbeta}
p_\beta(z,t) = p_0(z,t) + \beta \, p_1(z,t) +o(\beta)
\end{equation} 
At leading order, we find the following ordinary differential equation for $p_0$
\begin{equation}\label{adiabatic.approx}
\frac{d}{dz} \left( (z^2+\ell- t ) \, p_0(z,t)+ \frac{1}{2} \frac{d}{d z} p_0(z,t)\right) =0\,.
\end{equation}  
Equation \eqref{adiabatic.approx} is the same as Eq. \eqref{fokker.planck.statio} of the previous section 
with the parameter $a$ replaced by $-\ell +t$. The solution is
\begin{equation*}
p_0(z,t)  = 2\, J_0\left(\ell-  t\right) \, \int_{-\infty}^z du \, e^{2(\ell- t) (u-z) + \frac{2}{3} (u^3-z^3)}
\end{equation*} 
where $J_0$ is the stationary flux given by Eq. \eqref{flux.statio.J0} such that $\int_\R p_0(z,t) dz = 1$.  
Hence, recalling that $p_\beta(z,t)$ is the law of $Z_{4t/\beta}$, we conclude that the flux in $-\infty$ at time $t$ 
is $J_0\left(\ell-  \beta t/4\right)$ and thus,  
\begin{equation}\label{eq.Nbeta.expansion}
\rho_\beta(]-\infty;\ell])  = \int_0^{+\infty} J_\beta\left(-\infty,t \right) \, dt = \frac{4}{\beta} 
\int_{-\infty}^{\ell} J_0\left(u\right) \, du + O(1) \,, 
\end{equation}
where $O(1)$ is a correction of order $1$ when $\beta\to 0$. 
It is now straightforward to deduce the empirical eigenvalue density 
\begin{equation}\label{eigenvalue.density}
\rho_\beta(\ell) \sim_{\beta\to 0} \frac{4}{\beta}  J_0\left(\ell\right)\,. 
\end{equation}
This formula contrasts with the result found in \cite{lloyd,halperin,mckean} where the empirical eigenvalue density of the operator 
$\mathcal{L}_0$ restricted on a finite interval $[0;L]$ is proportional to the length $L$ of the time interval 
\begin{equation*}
n_L'(\ell)\sim_{L \to -\infty} L\, J_0'(\ell)\,. 
\end{equation*}
The computation leading to this result was in fact recalled above in \eqref{emp.eigenvalue.statio}. 

It is interesting to study the behaviors at $\ell \rightarrow \pm \infty$ 
of the empirical eigenvalue density $\rho_\beta(\ell)$ and of the integrated density $\rho_\beta(]-\infty;\ell])$.
This can be done from the integral forms with the saddle point route.  
When $\ell \rightarrow -\infty$, we have (see again end of Appendix \ref{moments.exit.time})
\begin{equation}\label{right.tail.intJ0}
\rho_\beta(\ell) \sim \frac{4}{\pi\beta} \, |\ell|^{1/2} \, \exp\left(-\frac{8}{3}|\ell|^{3/2}\right) \quad {\rm and} \quad  \rho_\beta(]-\infty;\ell])   
\sim \frac{1}{\pi\beta} \exp\left(-\frac{8}{3}|\ell|^{3/2}\right)\,. 
\end{equation}
On the other side, when $\ell \rightarrow +\infty$, we obtain 
\begin{equation*}
\rho_\beta(\ell) \sim  \frac{1}{\pi\beta} \, \ell^{1/2} \quad {\rm and} \quad   \rho_\beta(]-\infty;\ell])  \sim \frac{8}{3\pi\beta} \, \ell^{3/2}\,. 
\end{equation*}

\begin{figure}[h!btp] 
    \center
    \includegraphics[scale=0.6]{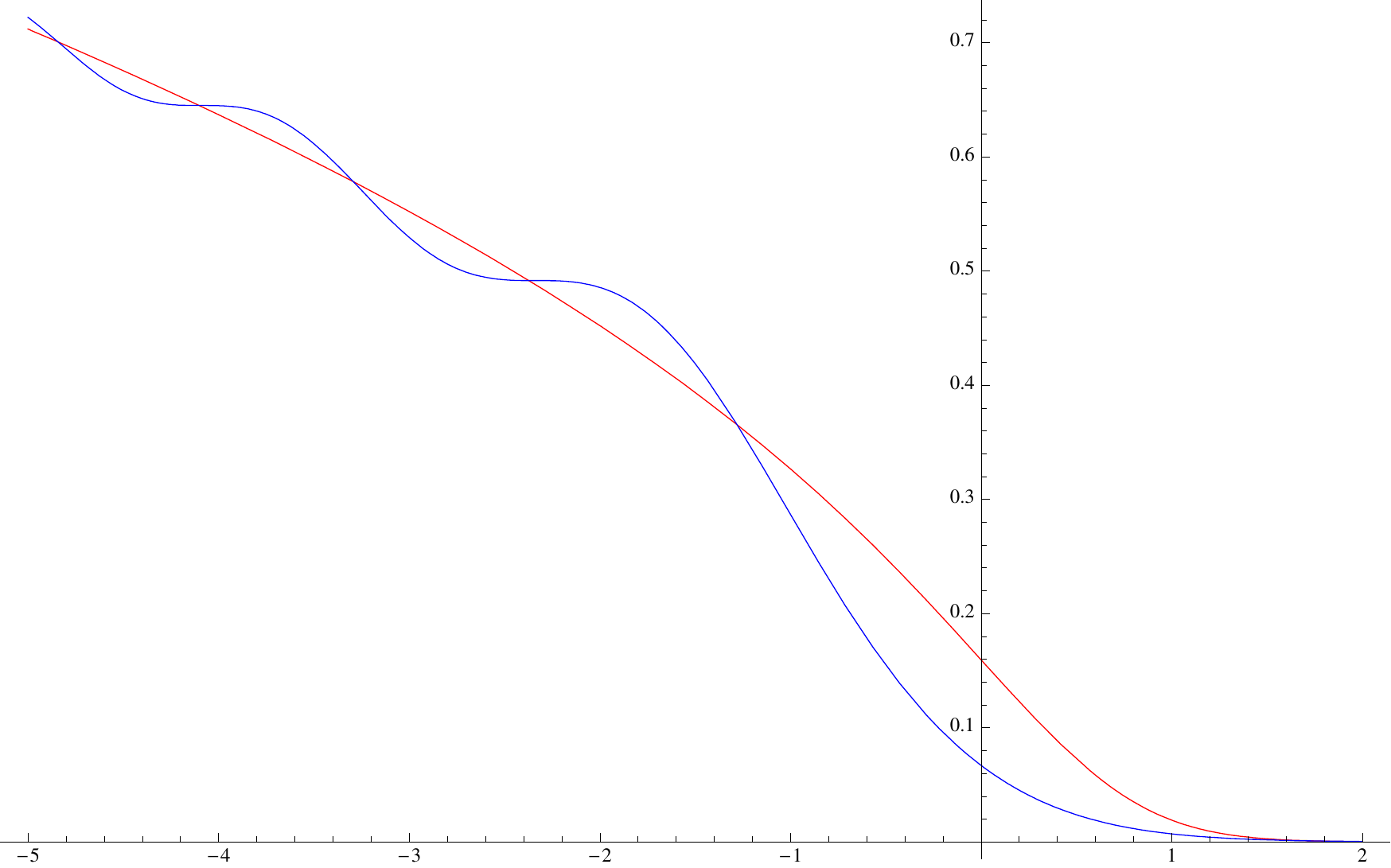}
    \caption{(Color online). Red non-oscillatory curve: $4 J_0(\lambda)$ as a function of $\lambda$. Blue oscillatory curve: 
     density $K_{{\rm Airy}}(\lambda,\lambda)= {\rm Ai}^\prime(\lambda)^2-\lambda {\rm Ai}^2(\lambda)$ at the edge of the spectrum for $\beta=2$ found by 
     Bowick and Br\'ezin. The bumps of the $\beta=2$ curve are reminiscent to the electrostatic repulsion between the particles.}\label{density.edge.scaling}
\end{figure}

We note the correct matching of the tails of the empirical eigenvalue density $\rho_\beta(\ell)$ as $\ell\to \pm \infty$ with those 
of the crossover density of eigenvalues (in random matrix models) near the edge scaling 
region of width $N^{-2/3}$ found by Bowick and Br\'ezin in \cite{bb} (see Figure \ref{density.edge.scaling}). 

We discuss the validity of perturbation theory as applied here in Appendix \ref{consistency.pt}.

\section{Applications to random matrix theory}  \label{rmt}
\subsection{Top eigenvalue of $\beta$-ensembles}
 In this subsection, we use the connection between random 
 matrices and stochastic linear operators to study the top 
 eigenvalue statistics of $\beta$-ensembles with an index $\beta_N$ depending of $N$ such that 
 \begin{equation}\label{hypbetaN}
 \beta_N\rightarrow_{N\to+\infty} \, 0 \, \quad {\rm and} \quad N \beta_N \rightarrow_{N\to+\infty} +\infty\,. 
 \end{equation}
 Example of such a sequence is $\beta_N=1/N^{\alpha}$ where $0<\alpha<1$. 
The tridiagonal random matrices introduced by Dumitriu and Edelman in \cite{dumitriu}, whose eigenvalues are distributed according to the jpdf 
$P_\beta$ defined in \eqref{eq.P_beta}, can be written as 
 \begin{align}\label{matrices.beta}
X_N =
\begin{pmatrix}
 \sqrt{2} \, g_1 & \chi_{(N-1)\beta_N} & \quad & \quad &  \\
 \chi_{(N-1)\beta_N} & \sqrt{2} \, g_2 & \chi_{(N-2)\beta_N} & \quad &  \\
 \quad  & \ddots & \ddots & \ddots &  \\
 \quad  & \quad & \chi_{2\beta} & \sqrt{2}\,  g_{N-1} & \chi_{\beta_N} \\
 \quad  & \quad & \quad & \chi_{\beta_N} & \sqrt{2}\, g_N
\end{pmatrix}
\end{align}
where the $g_k$ are independent Gaussian random variables with variance $1$ and where
the $\chi_{k\beta_N}$ are independent 
$\chi$ distributed random variables with $k\beta_N$ degrees of freedom and scale parameter $2$. 

In \cite{virag}, the authors prove that the largest eigenvalues of those tridiagonal random matrices 
converge in distribution to the low lying eigenvalues of the stochastic Airy operator $\mathcal{H}_\beta$ introduced above. 
Their proof is rather technical, but a simple heuristic of this convergence can be found in \cite{edelman}.
 
In the following, we apply this heuristic to analyse the present case with a vanishing repulsion 
coefficient $\beta_N$ satisfying \eqref{hypbetaN}. 
The main idea is that, under the assumption \eqref{hypbetaN}, the information on the top eigenvalues of the matrix $X_N$ is contained 
in the upper left sub matrices of $X_N$. In order to analyze the law of those eigenvalues, it is sufficient 
to consider only the entries with column or line numbers 
very small compared to the dimension $N$ of the matrix. 
Developing for $k\ll N$ the $\chi$ variables as $\chi_{\beta_N(N-k)} \approx \sqrt{\beta_N(N-k)} + \frac{h_k}{\sqrt{2}}$ where the $h_k$ are independent Gaussian variables, we obtain the following decomposition for the translated and rescaled matrix
%
\begin{align}\label{decomposition.matrix}
A_N&:= (N\beta_N)^{2/3} \left(2 I_N - \frac{X_N}{(N\beta_N)^{1/2}} \right)  = - (N\beta_N)^{2/3} \, \begin{pmatrix}
-2 & 1 & \quad & \quad \\
1 & -2 & 1 & \quad \\
\quad&\ddots&\ddots&\ddots
\end{pmatrix}\\ &+ \beta_N   \frac{1}{2} \frac{1}{(N\beta_N)^{1/3}}
\begin{pmatrix}
0 & 1 & \quad & \quad \\
1 & 0 & 2 & \quad \\
\quad&2 &\ddots&\ddots\\
\quad&\quad&\ddots&\quad
\end{pmatrix} \notag 
- \frac{1}{\sqrt{2}} \, (N\beta_N)^{1/6}
\begin{pmatrix}
2 g_1 & h_1 & \quad & \quad \\
h_1 & 2 g_2 & h_2 & \quad \\
\quad&\ddots&\ddots&\ddots
\end{pmatrix}+ o\left(\frac{1}{N}\right)\,.  
\end{align}
Using the same argument as in \cite{edelman},  
\eqref{decomposition.matrix} may be rewritten in short as 
\begin{equation}\label{eq.mat.op}
A_N \approx - \frac{d^2}{dt^2} + \beta_N \, t + 2\, b'_t\,. 
\end{equation}
In view of the previous section, it is important to keep the linear term $\beta_N t$ (even for $N\rightarrow +\infty$) 
so that the least eigenvalue of the operator 
on the right hand side of \eqref{eq.mat.op} remains finite. This term is small for small values of $t$ but gets large for
values of $t\gg \beta_N^{-1}$ and should not neglected.  

The top eigenvalue $\lambda_0^N$ of the matrix $X_N$ should therefore approximately satisfy, for large values of $N$, 
the equality in law 
\begin{align}\label{eq.top.eigenvalue}
(N\beta_N)^{2/3} \, \left(2 - \frac{\lambda_0^N}{(N\beta_N)^{1/2}}\right) \stackrel{(law)}{\approx} 4^{2/3}\, L_0^N
\end{align}
where $L_0^N$ is the least eigenvalue of the operator $-\frac{d^2}{dt^2}+\frac{\beta_N}{4} t + b_t'$. 
Using the previous results, Eq. \eqref{eq.top.eigenvalue} can be rewritten 
\begin{align}\label{final.result}
(N\beta_N)^{2/3} \left( \frac{\lambda_0^N}{(N\beta_N)^{1/2}}- 2  \right) \stackrel{(law)}{\approx} 
\left(\frac{3}{2}\right)^{2/3} \left(\ln \frac{1}{\pi\beta_N}\right)^{2/3}
+\left(\frac{2}{3}\right)^{1/3} \left(\ln \frac{1}{\pi\beta_N}\right)^{-1/3} \, G\,. 
\end{align}
where $G$ is a random variable distributed according to the Gumbel distribution. This result \eqref{final.result} 
was checked numerically with very good agreement.   

From this proposal that the tridiagonal random matrices in the $\beta$ ensembles may be regarded at the edge of the spectrum 
as finite different schemes of the stochastic Airy operator, we argue 
that the results proved in the continuous setting for the stochastic Airy operator 
can be extended to the discrete setting of random matrices. 
The convergence we obtained for the marginal distribution of any eigenvalue with fixed index $k$ in \eqref{other.vps.cv} 
can be rewritten in the context of $\beta_N$ ensembles as was done in \eqref{final.result} for the maximal eigenvalue. 
 
 At this point, it is tempting to conjecture, as in the continuous setting of the stochastic Airy operator, 
 that the maximal eigenvalues of $\beta_N$-ensembles have Poissonian statistics in the double scaling limit $N\to +\infty, 
 \beta_N\to 0$. Nevertheless the transition between Wigner/ Airy and Poisson statistics is not very well known up to now and 
 many questions have been asked in the literature of RMT (see e.g. \cite{cizeau,bordenave,victor,giulio} 
 for related questions on eigenvectors localization/delocalization).

\subsection{Break down at $\beta =2c/N$}
The discussion of the previous subsection breaks down if $\beta$ scales with $N$ as $\beta_N= 2c/N$ 
where $c$ is a positive constant. 
In such cases (and for even more rapidly decreasing $\beta_N$), the lower part of the matrix $X_N$ is no longer negligible 
compared to the upper part. One has to keep all the entries in the matrix $X_N$ and the approach through the 
stochastic linear operator is no longer valid.   

As mentioned in the introduction, the case $\beta=2c/N$ was studied before. 
In \cite{jp-alice}, the authors prove that, 
if $\beta$ scales with $N$ as $\beta=2c/N$ where $c$ is a positive constant, then the empirical eigenvalue density 
$1/N \sum_{i=1}^N \delta(\lambda-\lambda_i)$ of the matrix $X_N$ converges in the large $N$ limit to a continuous probability density  
$\rho_c$ with non compact support. The density $\rho_c$ is computed explicitly in \cite{jp-alice}
and one can easily recover the Gaussian tails of $\rho_c$ for $\lambda\to \pm\infty$  
\begin{align}\label{tail.rho.c}
\rho_c(\lambda) \sim C \lambda^{2c} e^{-\lambda^2/4}
\end{align} 
where $C$ is an explicit multiplicative constant. 

The reader may wonder whether Eq. \eqref{final.result} is coherent with the results found in \cite{jp-alice} where the limiting 
empirical eigenvalue density of the matrix $X_N$ was computed in the double scaling limit $\beta_N=2c/N$ with $N\to \infty$.  
 The developments of the previous subsection should a priori apply for any sequence $\beta_N$ such that $\beta_N\to 0$
and $N\beta_N\to +\infty$ as $N\to +\infty$. In particular, it should hold for $\beta_N= \ln N /N$. 
For such a $\beta_N$, Eq. \eqref{final.result} rewrites under the form 
\begin{align}\label{Gumbel.fluct}
\lambda_0^N = c_1 \sqrt{\ln N} + c_2 \frac{\ln \ln N}{\sqrt{\ln N}} + c_4 \frac{1}{\sqrt{\ln N}} +  c_3 \frac{1}{\sqrt{\ln N}} G\,,  
\end{align}
where $c_1,c_2,c_3,c_4$ are (explicit) constant (their values are irrelevant in the present discussion). 
This scaling form matches the one would  find for the maximum of $N$ particles allocated according to the 
density $\rho_c$ derived in \cite{jp-alice} when $\beta_N=2c/N$. 

We therefore conjecture that, even though $\lambda_0^N$ has Gumbel fluctuations for any sequences $\beta_N$ going to 
$0$ in the large $N$ limit, 
it has different scalings depending on whether 
\begin{itemize}
\item $\beta_N$ decreases faster than $\ln N/N$, then the centering and scaling 
in the convergence of $\lambda_0^N$ are as in \eqref{Gumbel.fluct}. 
\item $\beta_N$ decreases slower than $\ln N/N$, then the centering and scaling are as in Eq. \eqref{final.result}. 
\end{itemize}

\appendix
\section{Proofs of auxiliary results}\label{proof.theo}

{\it Proof of Proposition \ref{unicity.statio}.}

It is classical to show that the function $g_\alpha$ defined in Eq. \eqref{laplace.statio} is a weak solution to the 
boundary value problem Eq. \eqref{ode.g} and \eqref{boundary.g} \cite{ito}. 
Note simply that the differential equation  
is $\mathcal{G} g_\alpha = \alpha g_\alpha$ where $\mathcal{G}$ is the infinitesimal generator associated to the diffusion process 
$(Y_a(t))$. 
The boundary condition \eqref{boundary.g} is satisfied by $g_\alpha$ as defined in Eq. \eqref{laplace.statio}
since the diffusion, when starting from $y$, will explode in a time going to $0$ when $y\to -\infty$.

At this point, we do not know yet that the function $g_\alpha$ is a $C^2$ function. Nevertheless,  
we can (at least formally) 
perform the change of function $u(y) = \exp[-2(y^3/3-  ay)] g_\alpha'(y)$
and show that 
the function $g_\alpha$ satisfies 
the fixed point equation \eqref{pt.fixe.g}. 
This computation can be made rigorous by considering regularizations of $g_\alpha$ and taking the limit in the end.  
It immediately follows from Eq. \eqref{pt.fixe.g} that the function $g_\alpha$ is $C^{\infty}$. 

We can now prove that there exists a unique bounded $C^2$ function $g_\alpha$ satisfying the boundary 
value problem specified by Eq. \eqref{ode.g} and \eqref{boundary.g}. If $h$ is another solution, 
it is straightforward to check from It\^o's formula that the process $e^{-\alpha t} h(Y_a(t))$ is a martingale. 
Besides, it is bounded and hence we can apply the stopping time Theorem which yields 
\begin{equation*}
h(y) = \E_y[e^{-\alpha \zeta} \, {\bf 1}_{\{\zeta<\infty\}} \, h(Y_\zeta)]  =  \E_y[e^{-\alpha \zeta}] = g_\alpha(y)\,,
\end{equation*}  
where we have used the boundary condition at $y\to-\infty$ satisfied by $h$ as well as the fact that $\zeta <\infty$ almost surely.
The unicity is proved. 
\qed

{\it Proof of Theorem \ref{conv.exp}.}

First note that $g_{\alpha/m(a)}(y)$ is non increasing with respect to $y$. In addition we have $0\leq g_{\alpha/m(a)}(y) \leq 1$
and thus $ g_{\alpha/m(a)}(y)$ admits a limit when $y\to +\infty$. We will denote this limit $g_{\alpha/m(a)}(+\infty)$.

From Theorem \ref{unicity.statio}, $g_{\alpha/m(a)}(y)$ satisfies 
for any $\alpha\in (0;1)$ the fixed point equation 
\begin{equation}\label{pt.fixe.g.alphama}
g_{\alpha/m(a)}(y) = 1 - 2\frac{\alpha}{m(a)} \int_{-\infty}^{y} dx \int_x^{+\infty} du \exp\left(2a (x-u) + \frac{2}{3} (x^3-u^3)\right)\, g_{\alpha/m(a)}(u)\,. 
\end{equation}
Let us define recursively a sequence $R_n(y,a)$ such that $R_0(y,a)=1$ and for $n\geq 1$
\begin{equation}\label{int.Rn}
R_n(y,a) := \frac{2}{m(a)} \int_{-\infty}^{y} dx \int_x^{+\infty} du \exp\left(2a (x-u) + \frac{2}{3} (x^3-u^3)\right) R_{n-1}(u,a) \,. 
\end{equation}
To begin note that, for $\alpha \geq 0$, and for all $u \in\R$, 
\begin{equation}\label{ineq.g.0}
0 \leq g_{\alpha/m(a)}(u) = \E_u[e^{-\alpha \zeta/m(a)}] \leq 1
\end{equation}
Using the upper bound of Eq. \eqref{ineq.g.0}, we obtain from Eq. \eqref{pt.fixe.g.alphama} a new lower bound 
\begin{align}\label{ineq.g.1}
1-\alpha \, R_1(y,a) \leq g_{\alpha/m(a)}(y) \leq 1\,. 
\end{align} 
Using now the lower bound of Eq. \eqref{ineq.g.1}, we obtain from Eq. \eqref{pt.fixe.g.alphama} a new upper bound 
\begin{align*}
1-\alpha \, R_1(y,a) \leq g_{\alpha/m(a)}(y) \leq 1- \alpha\, R(y,a) + \alpha^2\, R_2(y,a) \,.
\end{align*}
We can check iteratively that for all $N$, 
\begin{align*}
\sum_{n=0}^{2N-1} (-1)^n \alpha^n R_n(y,a) \leq g_{\alpha/m(a)}(y) \leq \sum_{n=0}^{2N} (-1)^n \alpha^n R_n(y,a) \,. 
\end{align*}
In particular, we have by letting $N\to +\infty$ that, for any $\alpha \in (0;1)$, 
\begin{equation}\label{eqg.explicit}
g_{\alpha/m(a)}(y) = \sum_{n=0}^{+\infty} (-1)^n \alpha^n R_n(y,a)\,. 
\end{equation}
Note that \eqref{eqg.explicit} permits to deduce the values of all the moments $\E_y[\zeta^n]= n! \, m(a)^n \, R_n(y,a) $.

Now it suffices to prove that for all fixed $n\in \N$ and for any $y\in \R$, $R_n(y,a)$ converges to $1$ when $a\to +\infty$. This is the content of the following Lemma \ref{lemmeRn} whose proof can be found at the end of this section.
\begin{lemma}\label{lemmeRn}
Let $f(a)$ such that $a^{1/4}(f(a)+a^{1/2})\to +\infty$. 
Then, for each $n\in \N$ and for any $y\geq f(a)$, $R_n(y,a)$ converges to $1$ when $a\to +\infty$.    
\end{lemma}
It is then straightforward 
to see that 
\begin{equation*}
g_{\alpha/m(a)}(y,a) \rightarrow \frac{1}{1+\alpha}\,. 
\end{equation*}
This convergence holds for any $\alpha\in (0;1)$. 
The theorem follows since $1/(1+\alpha)$ is the Laplace transform of an exponential law with parameter $1$. 
\qed

\noindent {\it Proof of Lemma \ref{lemmeRn}.} 

First note that for each $n$, $R_n(y,a)$ is increasing with respect to $y$ and uniformly bounded by $1$. This can be seen 
easily by induction over $n$ recalling the expression of $m(a)(=m(a,+\infty))$ given in \eqref{eq.ma.explicit}. 
It follows that the $\lim_{y\to+\infty} R_n(y,a)$ exists and is denoted as usual 
$R_n(+\infty,a)$.   We proceed by induction over $n$ to show that $R_n(y,a)$ converges to $1$ for any $n$ and for any $y\geq f(a)$. 
By monotonicity, it is enough to prove that both $R_n(f(a),a)$ and $R_n(+\infty,a)$ 
converges to $1$ when $a\to +\infty$.

By the definitions of $R_1(y,a)$ (Eq. \eqref{int.Rn} for $n=1$) and of $m(a)$ (see Eq. \eqref{eq.ma.explicit}), 
it is clear that, for $y\in \R$ fixed and $a\to +\infty, R_1(y,a)\to 1$. 
Moreover, by monotonicity arguments, we easily check that this convergence holds uniformly on the interval $[f(a);+\infty[$.   

At step $n$, performing two changes of variables in the recursive integral expression Eq. \eqref{int.Rn} of $R_n(y,a)$,  
we obtain for all $a>0$, 
\begin{align*}
R_n(y,a) &= \frac{2}{m(a)}\, a \int_{-\infty}^{y/a^{1/2}} dx \int_{x}^{+\infty} du \exp\left(2a^{3/2} \left(\frac{1}{3}(x^3-u^3) - x + u \right) \right)\, 
R_{n-1}(a^{1/2} u,a)  \\
&= \frac{2}{m(a)} \frac{1}{a^{1/2}} \,  \exp\left(\frac{8}{3} a^{3/2}\right)  
\int\limits_{-\infty}^{a^{1/4} (y+a^{1/2})} dx \notag 
 \int\limits_{-2 a^{3/4}+x}^{+\infty} du \\&
 \exp\left( -2x^2 -2 u^2 +\frac{2}{3} (\frac{x^3}{a^{3/4}} - \frac{u^3}{a^{3/4}}) \right) R_{n-1}(a^{1/2} (1+\frac{u}{a^{3/4}}) ,a)\notag\,. 
\end{align*}
By monotonicity again, we simply need to consider $y\in \{f(a); +\infty\}$. 
From the second line and provided that $a^{1/4} (y+a^{1/2})\to +\infty$ (which is insured by assumption), 
we see that this integral is concentrated in the regions $x=O(1), u=O(1)$. 
One can in fact show that the other regions 
bring negligible contributions. 
In this regime $u=O(1)$ and we have, by the induction hypothesis, $R_{n-1}(a^{1/2} (1+\frac{u}{a^{3/4}}) ,a) \to 1$ uniformly for $u \in [-a^\epsilon;+\infty], \epsilon >0$ small.
The convergence of $R_n(f(a),a)$ and $R_n(+\infty,a)$ follow from the previous arguments and the asymptotic expression for 
$m(a)$ given in \eqref{asymp.ma}. 
\qed

\section{Proof of \eqref{nice.estimate} } \label{proof.paper.laure}

Setting $a=-\ell$ to work with a positive parameter and writing $Z_a$ instead of $Z_\ell$, 
recall the definition of the diffusion $Z_a$,
\begin{align*}
 dZ_a(t) = (\frac{\beta}{4} \, t + a -Z_a(t)^2)\, dt + dB(t)\,. 
\end{align*}
Recall also that $\P_{z,t}$ refers to the law of the diffusion $Z_a$ starting from $z$ at time $t$. We denote in this paragraph by $\zeta^{(u)}$ the first blowup time after time $u$ of $Z_a$.

Proposition 10 of \cite{laure} yields that 
\begin{lemma}\label{lemma.laure}
Fix $c' >0$. There exists a positive constant $c$ independent of $a$ and $\beta\in [0;1]$ such that for all $a$ large enough and for all $u \leq c' \ln(a)/(\beta \sqrt{a})$
\begin{align}\label{fine.estimate}
\P_{-\sqrt{a}+ \frac{(\ln a)^{1/4}}{a^{1/4}},u}\left[\zeta^{(u)} <\frac{\ln a}{\sqrt{a}}\right] \geq \exp(-c\sqrt{\ln a})\,. 
\end{align}

\end{lemma} 
\begin{rem}
 The constant in front of $\ln a/\sqrt{a}$ on the LHS is not optimal and can be replaced by any value strictly greater than $3/8$.
\end{rem}

\begin{proof}
For $x\in \R\cup \{-\infty\}$, set $T_x:= \inf\{t\geq 0: Z_a(t)=x\}$ and let $\delta:= (\ln a)^{1/4}/a^{1/4}$
and $\eps:= \sqrt{\ln a}/a^{1/4}$. Without loss of generality we can suppose that $B_0=0$. 

Moreover, to simplify the proof, we will restrict to the case $u=0$. The other cases follow the same lines 
(the only change to make is the value of $A$ introduced below \eqref{valuesAetB} which should be replaced by $A' := a + (1+c')\frac{\ln a}{\sqrt{a}}$).

Using the strong Markov property and the increasing property \cite{laure,virag}, we can lower bound the left hand side of \eqref{fine.estimate} by 
\begin{align}\label{three.proba}
\P_{-\sqrt{a}+ \delta,0}\left[ T_{-\sqrt{a}-\delta} < \frac{1}{\sqrt{a}} \wedge T_{-\sqrt{a}+2\delta} \right] \times 
\P_{-\sqrt{a}- \delta, \frac{1}{\sqrt{a}}}\left[ T_{-\sqrt{a}-\eps} < \frac{\ln a}{8 \sqrt{a}} \wedge T_{-\sqrt{a}}\right] \\
\times \P_{-\sqrt{a}- \eps, \frac{\ln a}{8\sqrt{a}}}\left[ T_{-\infty} < \frac{3}{8} \frac{\ln a}{\sqrt{a}} \right]\,. 
\end{align}

$\bullet$ The first probability gives the main cost. Using the comparison Theorem for sde (see 
\cite[Proposition 2.18]{karatzas} or \cite[Theorem (3.7), Chapter IX]{yor}), we see that the process $Z_a$ starting from 
$Z_a(0)= \sqrt{a}+\delta$ and until the stopping time $T_{-\sqrt{a}-2\delta}$ is stochastically dominated for $a$ large enough
 by the drifted Brownian motion $-\sqrt{a}+\delta + 4 \sqrt{a} \, \delta \, t + B(t)$.  
More precisely, this means that for all
$t \leq T_{-\sqrt{a}+2\delta}$, we have a.s. 
\begin{align*}
Z_a(t) \leq  -\sqrt{a}+\delta + 4 \sqrt{a} \, \delta \, t + B(t)\,. 
\end{align*}
Therefore, using the fact that $B(t) + 4 \sqrt{a} \, \delta \, t \geq B(t)$ for all $t\geq 0$, we can lower bound 
the first probability of \eqref{three.proba} by 
\begin{align}
\P\left[B(\frac{1}{\sqrt{a}}) < - 5 \frac{(\ln a)^{1/4}}{a^{1/4}} , \sup_{0\leq t \leq 1/\sqrt{a}} B(t) < \frac{(\ln a)^{1/4}}{a^{1/4}}\right]\,.
\end{align}
By the reflection principle \footnote{For $b \leq a$ and $t\geq 0$, 
$\P[ \sup_{0\leq s \leq t} B(s) \geq a, B(t) \leq b,] = \P[B_t \geq 2a-b]$.} 
and Brownian scaling, this later probability is 
\begin{align*}
\P[B(1)\leq - 5 (\ln a)^{1/4} ] 
-\P[ B(1) \geq 7 (\ln a)^{1/4}]&= \P[-7 (\ln a)^{1/4} \leq B(1) \leq -5 (\ln a)^{1/4} ] \\ 
& \geq c' \, \exp(-\frac{25}{2}\, \sqrt{\ln a})\,, 
\end{align*}
where the last inequality holds for $a$ large enough and for a positive constant $c'$.   
  
$\bullet$ For the second part, using again the comparison theorem for sde \cite[Proposition 2.18]{karatzas}, 
we can see that the process $Z_a$ starting from $Z_a(\frac{1}{\sqrt{a}})= -\sqrt{a}-\delta$ is almost surely 
below the Brownian motion for $t \in [\frac{1}{\sqrt{a}};  T_{-\sqrt{a}}]$. More precisely,   
\begin{align*}
Z_a(t) \leq - \sqrt{a} - \delta + B(t)- B(\frac{1}{\sqrt{a}})
\end{align*} 
for all $t\in  [\frac{1}{\sqrt{a}};  T_{-\sqrt{a}}]$.

Therefore the second probability is bounded from below by 
\begin{align*}
\P\left[B(\frac{\ln a}{8 \sqrt{a}}) \leq \delta - \eps, 
\sup_{0 \leq t \leq \frac{\ln a}{8\sqrt{a}}  } B(t) \leq \delta  \right]\,. 
\end{align*}
Again from reflexion principle and Brownian scaling property, this probability is equal to 
\begin{align*}
\P\left[-\sqrt{8}- \frac{\sqrt{8}}{(\ln a)^{1/4}} \leq B_1 \leq -\sqrt{8}+ \frac{\sqrt{8}}{(\ln a)^{1/4}} \right] =_{a\to +\infty}
O( \frac{1}{(\ln a)^{1/4}})\,. 
\end{align*} 
This handles the second probability of \eqref{three.proba} which does not contribute to the main cost as it 
is much larger than the first probability. 

$\bullet$ For the last probability, the idea of \cite{laure} is to compare 
the diffusion $Z_a$ with the solution of an (random) autonomous 
differential equation which can be computed explicitly. 
Set 
\begin{align}\label{def.Ga}
G_a(t):= Z_a(t) - B(t)\,. 
\end{align} 
 Denote 
 \begin{align*}
M := \sup_{0\leq t \leq \frac{3}{8}\frac{\ln a}{\sqrt{a}}} |B(t)| \,. 
\end{align*}
For $t\leq T_{-\sqrt{a}}\wedge \frac{3}{8}\frac{\ln a}{\sqrt{a}}$, we have 
\begin{align}\label{ineqGaM}
G_a(t) \leq -\sqrt{a} + M\,. 
\end{align} 
The function $G_a(t)$ satisfies the following (random) first order differential equation for $t> \frac{\ln a}{\sqrt{a}}$
\begin{align}\label{eq.diffGa}
G_a'(t) = \frac{\beta}{4}\, t+ a - G_a(t)^2 \left(1+\frac{B(t)}{G_a(t)}\right)^2\,, \quad G_a(\frac{\ln a}{\sqrt{a}}) 
= -\sqrt{a} - \frac{\sqrt{\ln a}}{a^{1/4}} \,. 
\end{align}  
Now consider the event 
\begin{align*}
\{ M\leq\frac{1}{2} \frac{\sqrt{\ln a}}{a^{1/4}}\} \,,
\end{align*}
under which we are able control the solution $G_a$ and its explosion time. 
Note that this event happens with a positive probability $p$ independent of $a$. 
Moreover, under this event, we have $T_{-\sqrt{a}}\geq \frac{3}{8} \frac{\ln a}{\sqrt{a}}$ so that the inequality \eqref{ineqGaM}
is valid for any $t\leq \frac{3}{8} \frac{\ln a}{\sqrt{a}}$. 

Thus, under the event $\{ M\leq\frac{1}{2} \frac{\sqrt{\ln a}}{a^{1/4}}\}$ and for $t\leq \frac{3}{8}\frac{\ln a}{\sqrt{a}}$, we easily check that 
\begin{align*}
 \left(1+\frac{B(t)}{G_a(t)}\right)^2 \geq 1- 2 \frac{M}{\sqrt{a}-M } \geq 1- \frac{\sqrt{\ln a}}{a^{3/4}- \frac{\sqrt{\ln a}}{2}}\,. 
\end{align*}  
Introduce now the function $H_a(t)$ solution of the (random) autonomous first order differential equation  
\begin{align}\label{eq.Ha}
H_a'(t) = A - B\, H_a(t)^2 \,, \quad H_a(0)= -\sqrt{a} - \frac{\sqrt{\ln a}}{a^{1/4}}  
\end{align}
with 
\begin{align}\label{valuesAetB}
A:= a+  \frac{\ln a}{\sqrt{a}} \quad  \mbox{ and } \quad B:=  1- \frac{\sqrt{\ln a}}{a^{3/4}- \frac{\sqrt{\ln a}}{2}}\,. 
\end{align}
It is easy to check (with the argument used to prove the increasing property) 
that, under the event $\{ M\leq\frac{1}{2} \frac{\sqrt{\ln a}}{a^{1/4}}\}$, the function $H_a$ 
dominates the function $G_a$ 
 a.s. for all $t\leq \frac{3}{8}\frac{\ln a}{\sqrt{a}}$,
\begin{align}\label{compGa.Ha}
G_a(\frac{\ln a}{\sqrt{a}}+ t) \leq H_a(t)\,. 
\end{align} 
But Eq. \eqref{eq.Ha} can be integrated with respect to $t$ leading to 
\begin{equation*}
\frac{H_a(t)+\sqrt{\frac{A}{B}}}{H_a(t)-\sqrt{\frac{A}{B}}} = C \, e^{2\sqrt{AB} t}
\end{equation*} 
where $C$ can be computed explicitly from the initial condition in $t=0$. From this expression, 
we see that the function $H_a$ explodes after a finite time $\tau$ which can be 
computed explicitly and bounded as
\begin{align}
\tau = \frac{1}{2\sqrt{AB}}\, \ln \frac{1}{C} \leq \frac{3}{8} \frac{\ln a}{\sqrt{a}} \,. 
\end{align} 
Using \eqref{def.Ga} and \eqref{compGa.Ha}, we deduce that under the event $\{ M\leq\frac{1}{2} \frac{\sqrt{\ln a}}{a^{1/4}}\}$, 
the diffusion $Z_a$ explodes in a time smaller 
than $(1+\frac{3}{8})\frac{\ln a}{\sqrt{a}}$. The Lemma is proved. 
\end{proof}

\section{Asymptotic of main integrals}\label{moments.exit.time}

With two consecutive change of variables, we obtain for $a>0$,
\begin{align}
m(a) 
&= \sqrt{2\pi}\,  a^{1/4} \int_{0}^{+\infty} \frac{dx}{\sqrt{x}} \exp\left(2a^{3/2} (x-\frac{1}{12} x^3 )\right) \notag \\
&= \frac{\sqrt{2\pi}}{\sqrt{a}} \exp\left(\frac{8}{3}a^{3/2}\right) \, \int_{-2a^{3/4}}^{+\infty} \frac{dy}{\sqrt{2+\frac{y}{a^{3/4}} }} \exp\left(-y^2 - \frac{1}{6} \frac{y^3}{a^{3/4}}\right) \label{integral.ma}\\
&\sim_{a\to + \infty} \frac{\pi}{a^{1/2}} \exp\left(\frac{8}{3} a^{3/2}\right)\,. \notag
\end{align}

With a more precise analysis of the integral \eqref{integral.ma}, we can in fact check that 
\begin{equation}\label{asymp.ma}
m(a) = \frac{\pi}{\sqrt{a}}  \exp\left(\frac{8}{3} a^{3/2} \right)   \left(1+\frac{5}{48} \frac{1}{a^{3/2}} + o(\frac{1}{a^3}) 
\right)\,,
\end{equation}
Then it is plain to deduce that, with $J_0(a) = 1/m(a)$,
\begin{equation}\label{left.tail.J0}
J_0(a)= \frac{a^{1/2}}{\pi}  
\exp\left(-\frac{8}{3} a^{3/2} \right)   \left(1-\frac{5}{48} \frac{1}{a^{3/2}} + o(\frac{1}{a^3})\right) \,, 
\end{equation}
and
\begin{equation}\label{intJ0.dl}
\int_{-\infty}^{\ell} J_0(t)\, dt \sim_{\ell \rightarrow -\infty} 
\frac{1}{4\pi} \exp\left(-\frac{8}{3} |\ell|^{3/2}\right) \left(1-\frac{5}{48} \frac{1}{|\ell|^{3/2}} 
+ o(\frac{1}{|\ell|^3})\right)\,. 
\end{equation}
Differentiating $J_0(a)$ with respect to $a$, we obtain for $a>0$
\begin{align}
J_0'(a)&= - 2 \, \sqrt{2\pi} \, J_0(a)^2 \, \int_0^{+\infty} dv\, \sqrt{v}\,  \exp\left(2a\,v-\frac{1}{6}v^3\right)\notag \\
&=  - 2 \, \sqrt{2\pi} \, J_0(a)^2 \, a^{3/4}  \, \int_0^{+\infty} dx\, \sqrt{x}\, \exp\left(2a^{3/2} (x-\frac{1}{12} x^3 )\right)\notag \\
&= - 2 \, \sqrt{2\pi} \, J_0(a)^2 \,\int_{-2a^{3/4}}^{+\infty}  dy\, \sqrt{2+\frac{y}{a^{3/4}}} \, \exp\left(-y^2 - \frac{1}{6} \frac{y^3}{a^{3/4}}\right)
\notag\\
&\sim_{a\to+\infty} -\frac{4}{\pi} \, a\, \exp\left(- \frac{8}{3} a^{3/2}\right)\,. \label{J0.prime}
\end{align}

\section{Consistency of perturbation theory}\label{consistency.pt}
Perturbation theory was used 
to obtain approximation at leading order of the solution of the
partial differential equation Eq. \eqref{eq.FPns}. 

The number $\rho_\beta(]-\infty;\ell])$ of eigenvalues below the 
level $\ell$ was then computed from the perturbative solution 
of the pde Eq. \eqref{eq.FPns}. 

In this section, we check the validity of perturbation theory by verifying that the correction terms remain negligible compared 
to the leading solution obtained. 

\subsection{Correction to the leading order}

We compute the $O(1)$ correction 
to the leading order in the expansion Eq. \eqref{eq.Nbeta.expansion} of $\rho_\beta(]-\infty;\ell])$ and show that, for all values of $\ell$,
perturbation theory leads to a consistent expansion at small $\beta$ 
of $\rho_\beta(]-\infty;\ell])$ such that the $O(1)$ correction is much smaller than 
the leading term of order $1/\beta$. 
In particular this holds for $\ell:=\ell_\beta$ scaling with $\beta$ as
$\ell_\beta \sim -\ln(1/\beta)^{2/3}$ when $\beta\to 0$. This was not obvious in the first place and needed to be checked in particular
for $\ell_\beta$ in the scaling region of the minimal eigenvalues at small $\beta$.    

The $O(1)$ correction denoted $\Gamma_\beta^1(\ell)$ such that 
\begin{equation}\label{expansion.N.L.beta}
\rho_\beta(]-\infty;\ell]) = \frac{4}{\beta}\int_{-\infty}^\ell J_0(u)\, du + \Gamma_{\beta}^1(\ell) + o(\beta) \,. 
\end{equation}
can be obtained by computing first the linear correction $J_1(z,t)$
to the flux $J_\beta^q(z,4t/\beta)= J_\beta^p(z,t)=J_0(\ell-t)+ \beta J_1(z,t) + o(\beta)$. 
The flux $J_1$ is related to the function $p_1$ introduced in Eq. \eqref{eq.pt.pbeta} and satisfies 
$J_1(z,t)= (z^2+\ell-t) p_1(z,t) +\frac{1}{2} \frac{d}{dz} p_1(z,t)$. We thus need to compute $p_1$.
This can be done by identifying the linear terms in $\beta$ on both sides of 
equation \eqref{eq.fp.pbeta}, we obtain the following ordinary differential equation for $p_1(z,t)$
\begin{align}\label{eq.J1}
 \frac{d}{d z} \left[ (z^2+\ell-t) \, p_1(z,t) +\frac{1}{2} \frac{d}{dz} p_1(z,t) \right]  = \frac{d}{dt} p_0^{\ell-t}(z)\,. 
\end{align}
By a further integration with respect to $z$, Eq. \eqref{eq.J1} becomes
\begin{equation*}
J_1(z,t) = \frac{d}{dt} \int_{-\infty}^z p_0^{\ell-t}(u) \, du + j_1(t)
\end{equation*}
where $j_1(t)$ is a constant which does not depend on $z$. 
As mentioned above, we are interested only in  $\lim_{z\rightarrow-\infty} J_1(z,t)= j_1(t)$, 
which can be computed easily using the normalization constraint $\int_\R p_1(z,t)dz=0$ 
\footnote{$p_1$ is the coefficient of the linear correction in $\beta$ to the probability density $p_\beta$: It should not bring mass in the integral.}. 
We find 
\begin{align*}
&j_1(t)= 4J_0'(\ell-t)  J_0(\ell-t) \int_{-\infty}^{+\infty} 
dz \int_{-\infty}^z du \int_{-\infty}^u dv  \int_{-\infty}^v dw \,e^{2(\ell-t)(u-z+w-v)+\frac{2}{3}(u^3-z^3+w^3-v^3)}  \\ 
&- 8J_0(\ell-t)^2 \int_{-\infty}^{+\infty} dz \int_{-\infty}^z du \int_{-\infty}^u dv  \int_{-\infty}^v dw \, (v-w) \,  e^{2(\ell-t)(u-z+w-v)+\frac{2}{3}(u^3-z^3+w^3-v^3)}  \,. 
\end{align*}
We can finally deduce the correction $\Gamma_\beta^1(\ell)$ 
as a function of the flux $j_1(t)$
\begin{equation*}
\Gamma_\beta^1(\ell) = \int_0^{+\infty} j_1(t) \, dt\,,
\end{equation*}
leading us to
\begin{align}
&\Gamma_\beta^1(\ell) = 4 \int_{-\infty}^\ell d\lambda J_0'(\lambda) J_0(\lambda) \int_{-\infty}^{+\infty} 
dz \int_{-\infty}^z du \int_{-\infty}^u dv  \int_{-\infty}^v dw \,e^{2\lambda(u-z+w-v)+\frac{2}{3}(u^3-z^3+w^3-v^3)}  \label{eq.N1}\\ 
&- 8\int_{-\infty}^\ell d\lambda J_0(\lambda)^2 \int_{-\infty}^{+\infty} dz \int_{-\infty}^z du \int_{-\infty}^u dv  \int_{-\infty}^v dw (v-w) e^{2\lambda(u-z+w-v)+\frac{2}{3}(u^3-z^3+w^3-v^3)} \notag \,. 
\end{align}
A careful analysis of this integral (see Appendix \ref{appendixB}) permits to extract the asymptotic behavior 
of $\Gamma_\beta^1(\ell)$ when $\ell\to-\infty$
\begin{align}\label{tail.N1}
\Gamma_\beta^1(\ell) \, \sim_{\ell\to-\infty} \, - \frac{3}{\pi} \, \ln |\ell| \, \exp\left(-4 |\ell|^{3/2}\right)\,.
\end{align}
Using this estimate, we can check that the asymptotic range $\ell_\beta(x)$ for $x\sim O(1)$ 
(going to $-\infty$ when $\beta\to 0$) 
of the minimal eigenvalue $L_0^\beta$ satisfies 
\begin{equation*}
\Gamma_\beta^1(\ell_\beta(x)) \propto -  \ln\left(\ln \frac{1}{\beta}\right) \, \beta^{3/2} \ll 
 \Gamma_\beta(\ell_\beta(x)) =O(1) \,. 
\end{equation*}
Thus we see that  the correction $\Gamma_\beta^1(\ell)$ remains negligible compared to $\rho_\beta(]-\infty;\ell])$ in the region 
$\ell\sim \ell_\beta(x)$.

\subsection{Other integrals} \label{appendixB}
If $\lambda<0$, then 
\begin{align*}
I_\lambda:=&\int_{-\infty}^{+\infty} 
dz \int_{-\infty}^z du \int_{-\infty}^u dv  \int_{-\infty}^v dw \,\exp \left( 2\lambda(u-z+w-v)+\frac{2}{3}(u^3-z^3+w^3-v^3)\right) \\ 
=& \int_{-\infty}^{+\infty} 
dz \int_{-\infty}^z du \int_{-\infty}^u dv  \int_{-\infty}^v dw \, \exp \left( 2|\lambda|(z-u+v-w)+\frac{2}{3}(u^3-z^3+w^3-v^3)\right) \\ 
= & |\lambda|^2 \int_{-\infty}^{+\infty} 
dz \int_{-\infty}^z du \int_{-\infty}^u dv  \int_{-\infty}^v dw \exp\left( 2 |\lambda|^{3/2} \left( (z-u+v-w) + \frac{1}{3} 
(u^3-z^3+w^3-v^3) \right) \right)\,. 
\end{align*}
In the limit $\lambda \rightarrow -\infty$, we can determine the leading order of this integral with the saddle point method. 
We have to compute the maximum of the function 
\begin{align*}
f(z,u,v,w) = (z-u+v-w) + \frac{1}{3} (u^3-z^3+w^3-v^3)
\end{align*} 
on the domain of integration $D:=\{(z,u,v,w)\in \R^4, u \leq z , v\leq u, w\leq v\}$. 
It is easily seen that it is reached for $z=1$, $u=v$, $w=-1$ and is equal to 
$4/3$. Let us perform the following change of variables 
\begin{align*}
z = 1+ \frac{x_1}{|\lambda|^{3/4}}, \quad u = x_2, \quad v = x_2- \frac{x_3}{|\lambda|^{3/2}}, \quad w= -1+\frac{x_4}{|\lambda|^{3/4}}\,. 
\end{align*} 
Up to small corrections, when $\lambda\rightarrow - \infty$, the integral is equivalent to 
\begin{align*}
|\lambda|^{-1} \exp\left(\frac{4}{3}|\lambda|^{3/2}\right) \int\limits_{-\infty}^{+\infty} dx_1 e^{-x_1^2} \int\limits_{-\infty}^{+\infty} dx_4 e^{-x_4^2}   \int\limits_{0}^{+\infty} dx_3 
\int\limits_{-1+ \frac{x_3}{|\lambda|^{3/2}}+\frac{x_4}{|\lambda|^{3/4}}}^{1+\frac{x_1}{|\lambda|^{3/4}}} dx_2 
  \exp\left( (x_2^2-1)  x_3 - x_2 \frac{x_3^2}{|\lambda|^{3/2}} \right)\,.
\end{align*}

Two contributions have to be accounted for. 
The variables $x_1$ and $x_4$ being fixed of order $1$, we look separately at the integration over $x_3$ with $x_2\leq 0$
and then with $x_2\geq 0$. 
\begin{align*}
&\int_{0}^{+\infty} dx_3  \int\limits_{-1+\frac{x_3}{|\lambda|^{3/2}} + \frac{x_4}{|\lambda|^{3/4}}}^0 dx_2 \exp\left( (x_2^2-1)  x_3 - x_2 \frac{x_3^2}{|\lambda|^{3/2}} \right) \\
& \sim_{\lambda\to -\infty}
\int_{0}^{+\infty} dx_3 \int\limits_{-1+\frac{x_3}{|\lambda|^{3/2}} }^0 dx_2 \exp\left( - (1- x_2^2)  x_3 \right)
 = \int_{0}^1 dx_2 \int_{0}^{(1-x_2)|\lambda|^{3/2}} dx_3 \exp\left(-(1-x_2^2) x_3\right) \\
&=   \int_{0}^1  \frac{dx_2}{1-x_2^2} \left(1- \exp\left(-(1-x_2^2) (1-x_2)|\lambda|^{3/2} \right) \right) 
= \int_{0}^1  \frac{dt}{t(2-t)} \left(1-\exp(-t^2 (2+t) |\lambda|^{3/2}) \right)  \\
&= \int_{0}^{|\lambda|^{3/4}} \frac{dx}{x\left(2-\frac{x}{|\lambda|^{3/4}}\right)} \left(1- \exp\left( - x^2 (2+\frac{x}{|\lambda|^{3/4}}) \right) \right)
\sim_{\lambda\to - \infty} \frac{3}{8} \ln |\lambda|\,. 
\end{align*}
Integrating now with respect to $x_1$ and $x_4$ leads to a contribution to the total integral $I_\lambda$ of order 
\begin{align*}
\pi \, |\lambda|^{-1} \, \exp\left(\frac{4}{3}|\lambda|^{3/2}\right) \, \frac{3}{8} \ln |\lambda|\,. 
\end{align*}

The region $x_2\geq 0$ has to be treated differently: Here we have to keep the correction term $-x_2x_3^2/|\lambda|^{3/2}$ in the exponential 
which prevents the integral to be infinite. 
\begin{align*}
\int\limits_{0}^{+\infty} dx_3 
&\int\limits_{0}^{1+\frac{x_1}{|\lambda|^{3/4}}} dx_2 
  \exp\left( (x_2^2-1)  x_3 - x_2 \frac{x_3^2}{|\lambda|^{3/2}} \right) \sim_{\lambda\to - \infty} 
  \int\limits_{0}^{+\infty} dx_3 
\int\limits_{0}^{1} dx_2 
  \exp\left( (x_2^2-1)  x_3 - x_2 \frac{x_3^2}{|\lambda|^{3/2}} \right) \\ 
 &= \frac{1}{|\lambda|^{3/2}} \int_0^{+\infty} dx_3 \int_0^{|\lambda|^{3/2}} dt \exp\left(-\frac{t}{|\lambda|^{3/2}}(2-\frac{t}{|\lambda|^{3/2}}) x_3 - 
 (1-\frac{t}{|\lambda|^{3/2}}) \frac{x_3^2}{|\lambda|^{3/2}} \right) \\
 &\sim_{\lambda\to -\infty}   \frac{1}{|\lambda|^{3/2}}  \int_0^{+\infty} dx \int_0^{|\lambda|^{3/2}} dt \exp\left(-\frac{2tx}{|\lambda|^{3/2}} - \frac{x^2}{|\lambda|^{3/2}}\right)\\ 
 &= \frac{1}{|\lambda|^{3/4}} \int_0^{+\infty} dy \int_0^{|\lambda|^{3/2}} dt  \exp\left(-\frac{2t}{|\lambda|^{3/4}} y  - y^2\right) \\
 &=  \frac{1}{2} \int_0^{+\infty} \frac{dy}{y} e^{-y^2}\left(1-e^{-2y |\lambda|^{3/4}}\right)\\
 &\sim_{\lambda\to -\infty}  \frac{1}{2} \int_0^{1} \frac{dy}{y} e^{-y^2}\left(1-e^{-2y |\lambda|^{3/4}}\right) = 
 \frac{1}{2} \int_0^{|\lambda|^{3/4}} \frac{dt}{t} e^{-t^2/|\lambda|^{3/2}} \left(1-e^{-2t}\right) \\
 &\sim_{\lambda\to -\infty} \frac{3}{8} \ln |\lambda|\,. 
\end{align*}   
Hence the region $x_2\geq 0$ leads to the same contribution as the region $x_2\leq 0$. We finally have 
\begin{align}\label{I.lambda}
I_\lambda \sim_{\lambda\to -\infty} \frac{3\pi}{4}\, \frac{\ln |\lambda|}{|\lambda|} \,  \exp\left(\frac{4}{3}|\lambda|^{3/2}\right) \,.
\end{align} 

Next we need to analyse the integral 
\begin{align*}
&K_\lambda:= \int_{-\infty}^{+\infty} dz \int_{-\infty}^z du \int_{-\infty}^u dv  
\int_{-\infty}^v dw (v-w) e^{2\lambda(u-z+w-v)+\frac{2}{3}(u^3-z^3+w^3-v^3)}\\
&= |\lambda|^{5/2} \int_{-\infty}^z du \int_{-\infty}^u dv  \int_{-\infty}^v dw (v-w) \exp\left( 2 |\lambda|^{3/2} \left( (z-u+v-w) + \frac{1}{3} 
(u^3-z^3+w^3-v^3) \right) \right)
\end{align*}
Up to small corrections, the integral $K_\lambda$ takes the form
\begin{align*}
&|\lambda|^{-1/2} \exp\left(\frac{4}{3}|\lambda|^{3/2}\right) \\& \int\limits_{-\infty}^{+\infty} dx_1 e^{-x_1^2} \int\limits_{-\infty}^{+\infty} dx_4 e^{-x_4^2}   \int\limits_{0}^{+\infty} dx_3 
\int\limits_{-1+ \frac{x_3}{|\lambda|^{3/2}}+\frac{x_4}{|\lambda|^{3/4}}}^{1+\frac{x_1}{|\lambda|^{3/4}}} dx_2 
 (x_2+1) \exp\left( (x_2^2-1)  x_3 - x_2 \frac{x_3^2}{|\lambda|^{3/2}} \right)\,.
\end{align*}
The integral in the second line can be analysed following the previous lines. 
We check that the mass in the integral is carried this time exclusively by the region $x_2\geq 0$.  
We obtain finally 
\begin{align}\label{K.lambda}
K_\lambda \sim_{\lambda\to-\infty} \frac{3\pi}{4}\, \frac{\ln |\lambda|}{ \sqrt{|\lambda|}} \,  \exp\left(\frac{4}{3}|\lambda|^{3/2}\right) \,.
\end{align}

Gathering \eqref{left.tail.J0}, \eqref{J0.prime}, \eqref{I.lambda} and \eqref{K.lambda}, we can deduce from \eqref{eq.N1} that
\begin{align*}
\Gamma_\beta^1(\ell) &\sim_{\ell\to-\infty} -\frac{18}{\pi}  \int_{-\infty}^\ell d\lambda\, \sqrt{ |\lambda|} \, \ln |\lambda|\,
\exp\left(-4|\lambda|^{3/2}\right) \\
&\sim_{\ell\to-\infty}  - \frac{3}{\pi} \, \ln |\ell| \, \exp\left(-4 |\ell|^{3/2}\right)\,.
\end{align*}

\end{document}